\documentclass[10pt, reqno]{amsart}
\usepackage{amscd,amssymb,epsfig}

\newtheorem{theorem}{Theorem}[section]
\newtheorem{lemma}[theorem]{Lemma}
\newtheorem{proposition}[theorem]{Proposition}
\newtheorem{corollary}[theorem]{Corollary}
\newtheorem*{thA}{Theorem 8.1}

\theoremstyle{definition}
\newtheorem{definition}[theorem]{Definition}

\theoremstyle{remark}

\numberwithin{figure}{section}
\numberwithin{table}{section}

\makeatletter

\@addtoreset{equation}{section}
\makeatother

\newcommand\abel{^{\operatorname{abel}}}
\newcommand\An{{\operatorname{Aut}(F_n)}}

\newcommand\Aut{{\operatorname{Aut}}}
\newcommand\bb{{\varsigma}}
\newcommand\bC{{\mathbb C}}
\newcommand\bZ{{\mathbb Z}}

\newcommand\bR{{\mathbb R}}
\newcommand\Coker{{\mbox{Coker}}}
\newcommand\dc{\overline{\partial}}
\newcommand\dd{{\partial}}
\newcommand\ddt{{\frac{d}{dt}\Bigr\vert_{t=0}}}

\newcommand\End{{\mbox{End}}}
\newcommand\family{{{\mathbb C}_{g}}}
\newcommand\fiber{{\int_{\text{fiber}}}}
\newcommand\Fn{{F_n}}
\newcommand\GL{{\mbox{GL}}}

\newcommand\harmonic{\mathcal{H}}
\newcommand\Hom{{\mbox{Hom}}}

\newcommand\IA{{\mbox{IA}}}

\newcommand\inter{int}
\newcommand\inv{^{-1}}

\newcommand\Ker{{\mbox{Ker}}}
\newcommand\La{{\Lambda}}
\newcommand\Lie{{\mbox{Lie}}}

\newcommand\Mgone{{\mathcal{M}_{g, 1}}}

\newcommand\modg{{{\mathbb M}_{g}}}
\newcommand\moduli{{{\mathbb M}_{g, 1}}}
\newcommand\mudot{{\stackrel{\centerdot}{\mu}}}

\newcommand\Mg{{\mathcal{M}_{g}}}
\newcommand\Mgstar{{\mathcal{M}_{g, *}}}

\newcommand\omegadot{{\stackrel{\centerdot}{\omega}}}
\newcommand\omegaone{{\omega_{(1)}}}

\newcommand\Sg{{\Sigma_g}}
\newcommand\stardot{{\stackrel{\centerdot}{\ast}}}

\newcommand\T{{\widehat{T}}}

\newcommand\thetadot{{\stackrel{\centerdot}{\theta}}}
\newcommand\wotimes{{\widehat{\otimes}}}
\newcommand\zbar{{\overline{z}}}

\begin{document}

\title[Harmonic Magnus Expansion]
{Harmonic Magnus Expansion 
on the Universal Family of Riemann Surfaces}
\author{Nariya Kawazumi}
\dedicatory{Dedicated to Professor Yukio Matsumoto on his sixtieth
birthday}
\address{Department of Mathematical Sciences\\
University of Tokyo \\Komaba, Tokyo 153-8914\\
Japan}
\email{kawazumi{\char'100}ms.u-tokyo.ac.jp}

\begin{abstract}
Let ${\mathbb M}_{g, 1}$, $g \geq 1$, 
be the moduli space of triples $(C, P_0, v)$ of  genus $g$, where  
$C$ is a compact Riemann surface of genus $g$, $P_0 \in C$, and 
$v \in T_{P_0}C\setminus\{0\}$. 
Using Chen's iterated integrals 
we introduce a higher analogue of the period matrix 
for a triple $(C, P_0, v)$,  {\it the harmonic Magnus expansion}. 
It induces a flat connection 
on a vector bundle over the space ${\mathbb M}_{g, 1}$, 
whose holonomy gives all the higher Johnson homomorphisms
of the mapping class group. 
The connection form, which is computed as an explicit 
quadratic differential, induces 
``canonical" differential forms representing
(twisted) Morita-Mumford classes and their higher relators 
on ${\mathbb M}_{g, 1}$.
In particular, we construct a family of twisted differential forms 
on ${\mathbb M}_{g, 1}$ representing 
the $(0, p+2)$-twisted Morita-Mumford class $m_{0, p+2}$ 
combinatorially parametrized by the Stasheff associahedron $K_{p+1}$.
\end{abstract}

\maketitle


\begin{center}
 Introduction
\end{center}

Let $\modg$ be the moduli space of compact Riemann surfaces 
of genus $g \geq 2$. The purpose of this and succeeding papers 
is to construct and study ``canonical" differential forms 
representing the Morita-Mumford classes (or the tautological 
classes) $e_i = (-1)^{i+1}\kappa_i$, $i \geq 1$, \cite{Mu} \cite{Mo1} 
on the moduli space $\modg$. 
For this purpose, in the present paper, we introduce 
a higher analogue of the period matrices of compact 
Riemann surfaces, the harmonic Magnus expansion, 
using Chen's iterated integrals \cite{C}. \par

It is well known the cohomology algebra $H^*(\modg; \bR)$ 
is naturally isomorphic to the algebra of all the real 
characteristic classes for oriented fiber bundles with 
fiber $\Sg$, a $2$-dimensional oriented connected closed 
$C^\infty$ manifold of genus $g$. 
This is the reason why the moduli space $\modg$ plays an 
important role also in differential topology. 
As was shown by Madsen and Weiss \cite{MW}, the stable part 
$* < g/3$ of the algebra $H^*(\modg; \bR)$ \cite{Har} is generated by the
Morita-Mumford classes 
$e_i$'s.\par

There are two classical approaches to constructing differential 
forms representing the classes $e_i$'s. Uniformization Theorem 
tells us the relative tangent bundle $T_{\mathbb{C}_g/\mathbb{M}_g}$ 
of the universal family $\pi: \family \to \modg$ 
of compact Riemann surfaces
has a canonical Hermitian metric, the hyperbolic metric. 
A notable work of Wolpert \cite{W}  
gives an explicit description of the differential forms on $\mathbb{M}_g$ 
representing the Morita-Mumford classes induced by the hyperbolic 
metric in terms of the resolvent of the hyperbolic Laplacian. 
On the other hand, from Grothendieck-Riemann-Roch Formula, 
the pullbacks of the Chern forms  on the Siegel upper half space 
$\frak H_g$ by the period matrix map represent all the {\bf odd}
Morita-Mumford classes. The first variation of the period matrices 
are given by Rauch's variational formula \cite{R} in terms of explicit 
quadratic differentials. \par

In order to construct differential forms representing 
{\bf all} the Morita-Mumford classes with no use of  
the hyperbolic metric, we introduce a higher analogue 
of the period matrices of compact Riemann surfaces. 
For simplicity, we consider the moduli space $\moduli$ of 
{\bf triples} $(C, P_0, v)$ of genus $g$ instead of the 
space $\mathbb{M}_g$. Here $C$ is a compact Riemann surface 
of genus $g$, $P_0 \in C$, and $v$ a non-zero tangent 
vector of $C$ at $P_0$. The space $\moduli$ is an aspherical  
$(3g-1)$-dimensional complex analytic manifold, and 
the fundamental group is equal to the mapping class group 
$\Mgone := \pi_0\operatorname{Diff}_+(\Sigma_g, p_0, v_0)$, 
where $p_0 \in \Sigma_g$, and $v_0 \in 
(T_\bR\Sigma_g)_{p_0}\setminus \{0\}$. 
The universal covering space is just the Teichm\"uller 
space $\mathcal{T}_{g, 1}$ for the topological triple 
$(\Sigma_g, p_0, v_0)$. \par

For any triple $(C, P_0, v)$ one can define the fundamental 
group of the complement $C\setminus \{P_0\}$ with the tangential 
basepoint $v$, which we denote by $\pi_1(C, P_0, v)$. 
If we choose a symplectic generator of $\pi_1(C, P_0, v)$, 
we can identify it with a free group of rank $2g$, $F_{2g}$. 
This induces a homomorphism $\Mgone \to \Aut(F_{2g})$, 
which is known to be an injection from a theorem of Nielsen. \par

Our construction is based on the following two facts: 
(1) We obtain the $i$-th  
Morita-Mumford class $e_i \in H^{2i}(\Mgone; \bR)$, 
$i \geq 1$, by contracting the coefficients 
of the $(0, 2i+2)$-twisted Morita-Mumford class $m_{0, 2i+2}
\in H^{2i}(\Mgone; H^{\otimes (2i+2)})$ by using the intersection 
product on the first homology group of the surface $\Sg$, 
$H := H_1(\Sg; \bR)$ \cite{Mo3} \cite{KM1} \cite{KM2}. 
(2) The $(0, p)$-twisted Morita-Mumford class $m_{0, p}$, $p \geq 3$,  
is an algebraic combination of the $p$-th power of 
the first Johnson map $[\tau^\theta_1]$, or equivalently, the $(0,
3)$-twisted  Morita-Mumford class $m_{0, 3}$ \cite{KM1} \cite{KM2}
\cite{Kaw2}.  From these facts, if one obtains a twisted 
differential form representing the class $m_{0, 3}$, 
then one can construct differential forms 
representing all the Morita-Mumford classes. 
In this paper we confine ourselves to studying the {\bf twisted} 
Morita-Mumford classes. 
Here we should remark a way to combine 
the classes $m_{0, 3}$'s to get the class $m_{0, p}$ 
is {\bf not} unique. There are many nontrivial 
relations among them, which come from the higher 
Johnson maps $\tau^\theta_s$, $s \geq 2$. 
In this paper we will show these relations are 
controlled by the Stasheff associahedron $K_{p+1}$ \cite{S}. 
For any $p \geq 1$, we construct a $p$-cocycle 
\begin{equation}
\theta^*Y_p \in C^*(K_{p+1}; \Omega^*(\moduli; H^{\otimes
(p+2)})) 
\end{equation}
of the cellular cochain complex of $K_{p+1}$ 
with values in the twisted de Rham complex of the moduli 
space $\moduli$. Here we regard the $\Mgone$-module $H^{\otimes (p+2)}$ 
as a flat vector bundle over the space $\moduli$ in an obvious way. 
The $p$-cocycle $\theta^*Y_p$ represents
$\frac{1}{(p+2)!}(-1)^{\frac12p(p+1)} m_{0, p+2}$
under the natural isomorphisms $H^p(C^*(K_{p+1}; \Omega^*(\moduli;\\
 H^{\otimes (p+2)}))) = H^p(\moduli; H^{\otimes(p+2)}) 
= H^p(\Mgone; H^{\otimes(p+2)})$, so that 
it can be interpreted as a 
canonical combinatorial family of differential forms 
representing the $(0, p+2)$-twisted Morita-Mumford class 
$m_{0, p+2}$. \par

The Johnson homomorphisms $\tau_p$, $p \geq 1$, are homomorphisms 
defined on a decreasing filtration $\{\mathcal{M}(p)\}^\infty_{p=0}$ 
of subgroups of the mapping class group introduced by Johnson \cite{J}. 
Independently of Johnson's work \cite{J} Harris \cite{H} defined 
the harmonic volume of a compact Riemann surface by using 
Chen's iterated integrals of harmonic forms and their relatives. 
Pulte \cite{P} proved a theorem of Torelli type 
for the pointed version of the harmonic volumes. 
Hain \cite{H} integrated these studies with Hodge theory 
to give an infinitesimal presentation of the Torelli group 
$\mathcal{M}(1)$. Kaenders \cite{Kae} proves a theorem of Torelli type 
for the Hodge structures of the fundamental group of 
 twice pointed Riemann surfaces. Our construction is more direct and 
elemantary than these Hodge-theoretic approaches. 
The Magnus expansion we will define in this paper 
is a higher analogue of Harris' harmonic volumes. 
This is the reason why we named it the {\bf harmonic} 
Magnus expansion. 
\par

In (0.1) 
we have two ingredients in the $p$-cocycle $\theta^*Y_p$, 
a cocycle $Y_p$ and a map $\theta$. 
Part I deals with the cocycle $Y_p$, which comes from 
some geometry about the automorphism group of a free group. 
Let $n \geq 2$ be an integer, and $\Fn$ a free group of rank $n$. 
We denote by $H := H_1(\Fn; \bR)$ the first real  
homology group of the group $\Fn$, and by $H^* := 
\Hom(H, \bR)$ its dual. The $p$-cocycle $Y_p$ inhabits 
the set $\Theta_n$ of all the real-valued Magnus expansions 
of the free group $\Fn$ in a generalized sense \cite{Kaw2}. 
A certain Lie group $\IA(\T)$ acts on the set $\Theta_n$ 
in a free and transitive way. This makes the set $\Theta_n$ 
a real analytic manifold, and induces the Maurer-Cartan form, 
which is a series of $1$-forms $\eta_p \in \left(\Omega^1(\Theta_n)
\otimes H^*\otimes H^{\otimes (p+1)}\right)^{\Aut(\Fn)}$, 
$p \geq 1$. The Johnson maps $\tau^\theta_p: \Aut(\Fn) \to 
H^*\otimes H^{\otimes (p+1)}$, $p \geq 1$, $\theta \in \Theta_n$, 
are given by Chen's iterated integrals of the forms $\eta_s$'s, 
$1 \leq s \leq p$. In other words, the Maurer-Cartan forms 
$\eta_s$'s define a flat connection on $\prod H^*\otimes H^{\otimes 
(p+1)}$ whose holonomy gives all the higher Johnson homomorphisms. 
Assembling the forms $\eta_s$'s we obtain 
the desired $p$-cocycle  
$Y_p \in C^*(K_{p+1}; \left(\Omega^*(\Theta_n)\otimes H^*\otimes
H^{\otimes (p+1)}\right)^{\Aut(\Fn)})$. \par

In Part II we construct an $\Mgone$-equivariant map 
$\theta: \mathcal{T}_{g, 1} \to \Theta_{2g}$, which we call 
{\bf the harmonic Magnus expansion}. For any triple 
$(C, P_0, v)$ of genus $g$, we obtain a canonical Magnus
expansion $\theta^{(C, P_0, v)}$ of the free group $\pi_1(C, P_0, v)$. 
Through the Poincar\'e duality we identify $H^* = H^1(C; \bR) 
= H^1(\pi_1(C, P_0, v); \bR)$ with $H = H_1(C; \bR) 
= H_1(\pi_1(C, P_0, v); \bR)$. In \S5, we construct 
a series of $1$-currents $\omega_{(p)}$ with coefficients in 
$H^{\otimes p}$, $p \geq 1$, such that $\omegaone$ comes from 
the harmonic forms on $C$ and $\omega = \sum\omega_{(p)}$ 
satisfies the integrability condition
$$
d\omega = \omega\wedge\omega - I\delta_0.
$$
Here $\delta_0$ is the delta current at $P_0$ and $I \in H^{\otimes 2}$ 
the intersection form of the surface $C$. Moreover we choose the 
current $\omega_{(p)}$ orthogonal to all the $d$-closed $1$-forms 
for any $p \geq 2$. 
The Magnus expansion $\theta^{(C, P_0, v)}$ is defined to be 
Chen's iterated integral of the connection form $\omega$ in \S 6, 
and gives the harmonic Magnus expansion $\theta: \mathcal{T}_{g, 1} \to
\Theta_{2g}$. In \S\S7 and 8 we compute the pullback $\theta^*\eta_p$, 
$p \geq 1$, of the Maurer-Cartan forms in terms of explicit quadratic
differentials.  Let $N: H^{\otimes m} \to H^{\otimes m}$
be a  linear map defined by $
N\vert_{H^{\otimes m}} :=\sum^{m-1}_{k=0}\begin{pmatrix}
1& 2& \cdots & m-1 & m\\
2& 3& \cdots & m & 1
\end{pmatrix}^k
$, and $\omega'$ the $(1, 0)$-part of the connection form $\omega$. 
Then we prove in Lemma 7.1 $N(\omega'\omega')$ is a meromorphic quadratic
differential  with coefficients in $\prod_{m \geq 2}H^{\otimes(m+2)}$ defined 
on $C$ with a pole of order $\leq 2$ at $P_0$. We denote by 
$N(\omega'\omega')_{(m)}$ the component of $N(\omega'\omega')$ 
with coefficients in $H^{\otimes m}$. Then we have 

\begin{thA} 
$$
(\theta^*\eta_p)_{[C, P_0, v]} = N(\omega'\omega')_{(p+2)} + 
\overline{N(\omega'\omega')_{(p+2)}}
\in (T^*_\bR\moduli)_{[C, P_0, v]}\otimes H^{\otimes (p+2)}
$$
for any $[C, P_0, v] \in \moduli$ and $p \geq 1$ . 
Here we regard the quadratic differential $N(\omega'\omega')_{(p+2)}$ 
as a $(1, 0)$-cotangent vector at $[C, P_0, v] \in \moduli$ in a
natural way.
\end{thA}

The second term of $N(\omega'\omega')_{(2)} = 2\omegaone'\omegaone'$
is just the first variation of the  period matrices given by 
Rauch's variational formula \cite{R}. \par
In \S8 we show the $1$-form $\theta^*\eta_1$ has coefficients 
in $\La^3H \subset H^{\otimes 3}$. 
Let $U := \La^3H/H$ be the cokernel of the
injection
$\frak q^H: H \to \La^3H$, $Z \mapsto Z\wedge I$, and 
$\frak p^U: \La^3H \to U$ the natural projection. 
In Theorem 8.3 we prove $\eta^U_1 := \frak p^U\theta^*\eta_1$ 
can be regarded as a $1$-form on $\modg$ representing the 
extended Johnson homomorphism $\tilde k \in H^1(\modg; U)$ 
introduced by Morita \cite{Mo2}. The $1$-form $\eta^U_1$ is 
exactly the first variation of the harmonic volumes 
given by Harris \cite{H}. Hain and Reed \cite{HR} introduced the 
$1$-form $\eta^U_1$ in the context of Hodge theory, and 
studied a $2$-form representing the first Morita-Mumford class 
$e_1$ obtained from $\eta^U_1$. 
As was pointed out by Harris, the form $\eta^U_1$ 
vanishes along the hypereliptic locus $\mathcal{H}_g \subset \modg$. 
This implies all the differential 
forms representing the Morita-Mumford classes derived from 
 $\eta^U_1$ vanish along the locus $\mathcal{H}_g$ in contrast to 
the differential forms coming from 
the hyperbolic metric. 
The results in \S\S4-7 were announced in \cite{Kaw3}.

\par
\smallskip\noindent
{\it Acknowledgements}: 
The author would like to thank Shigeyuki Morita, Kiyoshi Igusa, 
Richard Hain, Kiyoshi Ohba, Hiroshi Ohta, Ryushi Goto, Masahiko
Yoshinaga, Yuuki Tadokoro, Takao Satoh, Yusuke Kuno and Masatoshi Sato 
for imspiring discussions. In particular, the author would like 
to express his sincere gratitude to Yukio Matsumoto, 
the author's former supervisor.
\par

\tableofcontents

\part{Geometry of Magnus Expansions}

\section{Iterated Integrals and Maurer-Cartan Forms}

We review iterated integrals introduced by K.-T. Chen
\cite{C}  and the Maurer-Cartan form on a Lie group. 
In \S 2 we will describe the Johnson map of the automorphism 
group of a free group as an iterated integral of 
a certain kind of Maurer-Cartan forms on the space 
of all the ($\bR$-valued) Magnus expansions of 
the free group.
Moreover, in Part II, the harmonic Magnus expansion will be 
defined as an iterated integral over the universal family 
of compact Riemann surfaces. 
\par

Our definition of iterated integrals is slightly different from 
Chen's original one \cite{C}. 
In fact, for a topological space $M$ and paths 
$\ell_0$, $\ell_1: [0, 1] \to M$ with $\ell_0(1) = \ell_1(0)$, 
we define the product $\ell_1\cdot\ell_0: [0, 1] \to M$ by 
\begin{equation}
(\ell_1\cdot\ell_0)(t) := 
\left\{ \begin{array}{ll}
\ell_0(2t), \quad& \text{if $0 \leq t \leq 1/2$,}\\
\ell_1(2t-1), \quad& \text{if $1/2 \leq t \leq 1$,}
\end{array}\right.
\end{equation}
which is in reverse order to \cite{C} 1.5, p.221. \par
Let $T$ be an associative $\bR$-algebra with a unit element $1$ 
with a decreasing filtration of two-sided ideals $T_p$, $p\geq 1$, 
of finite real codimension. 
We assume $T_p\cdot T_q \subset T_{p+q}$ for any $p$ and $q \geq 1$, 
and the natural map $T \to \varprojlim_{p\to \infty}T/T_p$ is 
an isomorphism
\begin{equation}
T \cong \varprojlim_{p\to \infty}T/T_p. 
\end{equation}
Let $f_i(t)$, $1 \leq i \leq m$, be a continuous function defined
on  a closed interval $[a, b]$, $a, b \in \bR$, with values in
the ideal 
$T_1$. 
We have a $T_1$-valued $1$-form $\varphi_i := f_i(t) dt$ on the 
interval $[a, b]$.  
The iterated integral of $\varphi_1, \varphi_2, \dots,
\varphi_m$  along $[a, b]$ is defined by 
$$
\int^b_a\varphi_1 \varphi_2 \cdots \varphi_m
:= \int_{b \geq t_1\geq \cdots \geq t_m \geq a} 
f_1(t_1)f_2(t_2)\cdots f_m(t_m)dt_1dt_2\cdots dt_m.
$$
If $m \geq 2$, we have
\begin{equation}
\int^b_a\varphi_1\varphi_2\cdots\varphi_m = 
\int^b_af_1(t_1)\left(\int^{t_1}_a\varphi_2\cdots\varphi_m
\right)dt_1. 
\end{equation}
Now suppose $f_1(t) = f_2(t) = \cdots = f_m(t) =: f(t)$, i.e.,
$\varphi_1 = \varphi_2 = \cdots = \varphi_m =: \varphi$. 
Then
$$
F(t) :=
1+{\sum}^\infty_{m=1}\int^t_a
\overbrace{\varphi\varphi\cdots\varphi}^{\mbox{\small $m$ times}}
$$
converges because of the assumption (1.2) and $f(t) \in T_1$. 
From (1.3) the function $F(t)$ is the unique solution of 
the initial value problem
\begin{equation}
\left\{ \begin{array}{ll}
& \displaystyle{\frac{d}{dt}}F(t) = f(t)F(t),\\
& F(a) = 1. \end{array}\right.
\end{equation}
\par

Let $M$ be a $C^\infty$ manifold, and $E$ a $C^\infty$ vector
bundle over $M$. 
We denote by $C^\infty(M; E)$ the space of all the $C^\infty$ 
sections of $E$ over $M$. 
If $E$ is a flat vector bundle over $M$, then we denote 
the twisted de Rham complex with coefficients in $E$ by 
$$
\Omega^*(M; E) = {\bigoplus}^\infty_{q=0}\Omega^q(M; E) =
{\bigoplus}^\infty_{q=0}C^\infty(M; (\La^qT^*M)\otimes E),
$$
where $\La^qT^*M$ is the
$q$-cotangent bundle of $M$, $q \geq 0$. 
Moreover we write simply $\Omega^q(M)\wotimes T_p$ 
for $\Omega^q(M; M\times T_p)$, $p \geq 1$, which is equal to 
$\varprojlim_{m\to\infty}(\Omega^q(M)\otimes (T_p/T_m))$. 
 For any piecewise $C^\infty$ path
$\ell: [0, 1] \to M$ and  any $T_1$-valued $1$-forms $\varphi_1,
\varphi_2, \dots,
\varphi_m \in \Omega^1(M)\wotimes T_1$ the iterated 
integral is defined by
$$
\int_\ell\varphi_1\varphi_2\cdots\varphi_m = 
\int^1_0({\ell}^*\varphi_1)({\ell}^*\varphi_2)
\cdots({\ell}^*\varphi_m).
$$ 
Let $h \in \Omega^0(M)\wotimes T$ be a $T$-valued function 
with $dh \in \Omega^1(M)\wotimes T_1$.
The following formulae are well-known, and easy to prove.
\begin{eqnarray}
&& \int_\ell(dh)\varphi_1\varphi_2\cdots \varphi_m =
h(\ell(1))\int_\ell\varphi_1\varphi_2\cdots \varphi_m - \int_\ell
(h\varphi_1)\varphi_2\cdots \varphi_m \nonumber
\\
&& \int_\ell\varphi_1\cdots\varphi_{i-1}(dh)\varphi_{i+1}\cdots
\varphi_m \\
&=&
\int_\ell\varphi_1\cdots\varphi_{i-2}(\varphi_{i-1}h)\varphi_{i+1}\cdots
\varphi_m -
\int_\ell\varphi_1\cdots\varphi_{i-1}(h\varphi_{i+1})\varphi_{i+2}\cdots
\varphi_m
\nonumber\\
&& \int_\ell\varphi_1\varphi_2\cdots \varphi_m (dh) =
\int_\ell\varphi_1\varphi_2\cdots (\varphi_mh) - \left(\int_\ell
\varphi_1\varphi_2\cdots \varphi_m\right)h(\ell(0))
\nonumber
\end{eqnarray}
Any $T_1$-valued $1$-form $\varphi \in \Omega^1(M)
\wotimes T_1$ defines the connection
$$
\nabla: \Omega^0(M)\wotimes T 
\to \Omega^1(M)\wotimes T, \quad
u \mapsto du - \varphi u
$$
on the product bundle $M \times T$. 
For any piecewise $C^\infty$ path $\ell: [0, 1] \to M$ 
we denote $\ell_t(s) := \ell(ts)$, $0 \leq t, s\leq 1$. 
Then, from (1.4), the iterated integral
$$
F_\ell(t) = 1+{\sum}^\infty_{m=1}\int_{\ell_t}
\overbrace{\varphi\varphi\cdots\varphi}^{\mbox{$m$ times}}
$$ 
is the horizontal lift of $\ell$ through $F_\ell(0) = 1$. 
If the connection $\nabla$ is flat, i.e., $\varphi$ satisfies 
the integrability condition 
\begin{equation}
d\varphi = \varphi\wedge\varphi,
\end{equation}
then $F_\ell(1)$ is invariant under any homotopy of $\ell$
fixing the endpoints $\ell(0)$ and $\ell(1)$.\par

Next we recall the Maurer-Cartan formula. 
Let $G$ be a Lie subgroup of the general linear group 
$\GL(V)$ for a $\bR$-vector space $V$, and $\frak g$ 
its Lie algebra $\frak g := \Lie(G)$. 
Suppose the manifold $M$ has a left free transitive 
$C^\infty$ action of $G$. 
Then $M$ is diffeomorphic to $G$. The global 
trivialization
$$
M\times \frak g \overset\cong\to TM, \quad
(x, u) \mapsto \ddt\exp(tu)x
$$
can be regarded as a $\frak g$-valued $1$-form 
$\eta \in \Omega^1(M)\wotimes\frak g$, which is, 
by definition, the Maurer-Cartan form for the action 
of $G$ on $M$. \par

The Maurer-Cartan form $\eta$ satisfies the integrability 
condition (1.6). In fact, choose a point $x_0 \in M$, and 
define a $C^\infty$ map
$$
F = F^{x_0}: M \to G \subset \GL(V)\subset\operatorname{End}(V)
$$
by $F(gx_0) = g$ for $g \in G$. 
Then, for any $u \in \frak g$ and $g \in G$, we have
\begin{eqnarray*}
&&((dF)F\inv)\left(\ddt\exp(tu)gx_0\right) =
\left(\ddt\exp(tu)g\right)g\inv\\ 
&=&\,\,\ddt\exp(tu) = u =
\eta\left(\ddt\exp(tu)gx_0\right),
\end{eqnarray*}
that is,
\begin{equation}
dF = \eta F \in \Omega^1(M) \wotimes
\End(V).
\end{equation}
This implies $0 = ddF = d(\eta F) = (d\eta -
\eta\wedge\eta)F$. Hence we obtain
\begin{equation}
d\eta = \eta\wedge\eta 
\in \Omega^2(M)\wotimes
\End(V),
\end{equation}
which is just the Maurer-Cartan formula.\par
From (1.4) and (1.7) we have 
\begin{equation}
g = 1+{\sum}^\infty_{m=1}\int^{gx_0}_{x_0}
\overbrace{\eta\eta\cdots\eta}^{\mbox{$m$ times}}
 \in \End(V)
\end{equation}
for any $g \in G$, where the iterated integral on the 
right hand side is along any piecewise $G^\infty$ path 
from $x_0$ to $gx_0$.

\section{Magnus Expansions and Johnson Maps}

Following \cite{Kaw2} we recall the notion of 
Magnus expansions of a free group in a generalized sense. 
A certain Lie group $\IA(\T)$ acts on the space of 
Magnus expansions, which yields the Johnson maps 
defined on the whole automorphism group of the free group.
\par

Let $n \geq 2$ be an integer, $\Fn$ a free group of rank $n$
with free basis $x_1, x_2, \dots, x_n$
$$
\Fn = \langle x_1, x_2, \dots, x_n\rangle.
$$
We denote the first real homology group of $\Fn$ by 
$$
H := \Fn\abel\otimes_\bZ \bR = H_1(\Fn; \bR),
$$
$[\gamma] := (\gamma \bmod [\Fn, \Fn])\otimes_\bZ 1 \in H$ 
for $\gamma \in \Fn$, and $X_i := [x_i] \in H$ 
for $i$, $1 \leq i\leq n$. 
For the rest of the paper we write simply $\otimes$ and 
$\Hom$ for the tensor product 
and the homomorphisms over the real numbers $\bR$, 
respectively. 
The completed tensor algebra generated by $H$ 
$$
\T = \T(H) : = {\prod}^\infty_{m=0} H^{\otimes m}
$$
is equal to the ring of noncommutative formal power series 
$\bR\left<\left<X_1, X_2, \dots, X_n\right>\right>$. 
The two-sided ideals
$$
\T_p := {\prod}_{m\geq p} H^{\otimes m}, \quad p \geq 1,
$$
give a decreasing filtration of the algebra $\T$, 
which satisfies the condition (1.2). 
For each $m$ we regard $H^{\otimes m}$ as a subspace of $\T$ 
in an obvious way. So we can write 
$$
z = \sum^\infty_{m=0} z_m 
= z_0 + z_1 + z_2 + \cdots + z_m + \cdots
$$
for $z = (z_m) \in \T$, $z_m \in H^{\otimes m}$.
The subset $1+\T_1$ is 
a subgroup of the multiplicative group of the algebra $\T$. 

\begin{definition}[\cite{Kaw2} Definition 1.1] 
A map $\theta: \Fn \to 1 +\T_1$ is  a (real-valued) Magnus
expansion of the free group $\Fn$, if 
\begin{enumerate}
\item $\theta: \Fn \to 1 + \T_1$ is a group homomorphism, and
\item $\theta(\gamma) \equiv 1 + [\gamma] \pmod{\T_2}$ for any
$\gamma 
\in \Fn$.
\end{enumerate}
We denote by $\Theta_n = \Theta_{n, \bR}$ 
the set of all the (real-valued) Magnus expansions. 
\end{definition}
\par
Now we consider the group $\Aut(\T)$ of all the
filtration-preserving $\bR$-algebra automorphisms 
of the algebra $\T$. 
Here an $\bR$-algebra automorophism $U$ of $\T$ 
is defined to be {\bf filtration-preserving} 
if $U(\T_p) = \T_p$ for each $p \geq 1$. 
The group $\Aut(\T)$ is a (projective limit of)
Lie group(s) in a natural way. 
We denote by $\vert U\vert \in \GL(H)$ the 
isomorphism of $H = \T_1/\T_2$ induced by $U \in \Aut(\T)$. 
We define 
$$
\IA(\T) := \Ker \vert\cdot\vert \subset \Aut(\T),
$$
which is a closed subgroup of the Lie group $\Aut(\T)$. \par

If $U \in \IA(\T)$ and $\theta \in \Theta_n$, 
then the composite $U \circ \theta: \Fn\overset\theta\to
1+\T_1\overset{U}\to1+\T_1$ is also a Magnus expansion. 
This means the group $\IA(\T)$ acts on the set $\Theta_n$. 
The following was proved in \cite{Kaw2} Theorem 1.3 (2).

\begin{theorem} The action of the group $\IA(\T)$ 
on the set $\Theta_n$ is free and transitive.
\end{theorem}

It enables us to regard the set $\Theta_n$ as 
a (projective limit of) $C^\infty$ manifold(s). 
We denote the Maurer-Cartan form for the action by
$$
\eta \in \Omega^1(\Theta_n) \wotimes\,\Lie\,\IA(\T).
$$

In order to study the Lie algebra $\Lie\,\IA(\T)$ 
we introduce the $\bR$-algebra $\End(\T)$ consisting of all 
the $\bR$-linear maps $u: \T \to \T$ satisfying 
$u(\T_m) \subset \T_m$ for each $m\geq 1$. 
$\IA(\T)$ and $\Lie\,\IA(\T)$ can be regarded as subsets 
of $\End(\T)$ in a natural way. 
A linear map $u \in \End(\T)$ belongs to $\Lie\,\IA(\T)$, 
if and only if $u(\T_1) \subset \T_2$ and $u$ is a derivation, 
i.e., it satisfies
\begin{equation}
u(ab) = u(a)b + au(b)
\end{equation}
for any $a$ and $b \in \T$. 
Hence the restriction map to $H$
$$
\Lie\,\IA(\T) \to \Hom(H, \T_2), \quad
u \mapsto u\vert_H
$$
is an $\bR$-linear isomorphism. 
If we denote $H^* = \Hom(H, \bR)$, then we have an isomorphism
\begin{equation}
\Lie\,\IA(\T) \cong \Hom(H, \T_2) = H^*\otimes \T_2 
= {\prod}^\infty_{p=1}H^*\otimes H^{\otimes (p+1)}.
\end{equation}
Immediately from (2.1) follows

\begin{lemma} For $u$ and $v \in \Lie\,\IA(\T) \subset 
\End(\T)$ denote
$$
u\vert_H = {\sum}^\infty_{p=1} u_p, \quad
v\vert_H = {\sum}^\infty_{p=1} v_p, \quad\text{and}\quad
uv\vert_H = {\sum}^\infty_{p=1} w_p,
$$
where $u_p$, $v_p$ and $w_p \in H^*\otimes H^{\otimes (p+1)}$. 
Then we have 
$$
w_p = \sum^{p-1}_{s=1}(\underbrace
{u_s\otimes 1\otimes
\cdots\otimes 1 + 1 \otimes u_s\otimes \cdots \otimes 1 + \cdots
+  1\otimes \cdots\otimes 1\otimes u_s}_{p-s+1}) \circ v_{p-s}.
$$
\end{lemma}
Let $\eta_p \in C^{\infty}(\Theta_n; T^*\Theta_n)\otimes 
H^*\otimes H^{\otimes (p+1)}$ be the $p$-th component of 
the Maurer-Cartan form 
$$
\eta = \sum^\infty_{p=1}\eta_p \in \Omega^1(\Theta_n)\wotimes 
\Lie\,\IA(\T) = \prod^\infty_{p=1}\Omega^1(\Theta_n)
\otimes  H^*\otimes
H^{\otimes (p+1)}.
$$
By Lemma 2.3 the Maurer-Cartan formula (1.8) is 
equivalent to 
\begin{lemma}
\begin{eqnarray}
& d\eta_1 = 0, \nonumber\\
& d\eta_2 = (\eta_1\otimes 1 + 1\otimes \eta_1)\eta_1\nonumber\\
& d\eta_p = \sum^{p-1}_{s=1}(\underbrace
{\eta_s\otimes 1\otimes
\cdots\otimes 1 + \cdots
+  1\otimes \cdots\otimes 1\otimes \eta_s}_{p-s+1})\eta_{p-s}.
\end{eqnarray}
\end{lemma}
As will be discussed in \S3, the formula (2.3) suggests 
a close relation between the Stasheff associahedron and the
Maurer-Cartan form $\eta$.\par
\medskip
Now we consider the automorphism group of the group $\Fn$, 
$\Aut(\Fn)$. We denote by $\vert\varphi\vert \in \GL(H)$ 
the induced map on $H = H_1(\Fn; \bR)$ by $\varphi \in 
\Aut(\Fn)$. 
Fix a Magnus expansion $\theta \in \Theta_n$. 
Then {\bf the total Johnson map} induced by $\theta$
$$
\tau^\theta: \Aut(\Fn) \to \IA(\T), \quad
\varphi \mapsto \tau^\theta(\varphi)
$$
is defined by 
\begin{equation}
\tau^\theta(\varphi)\inv\circ\theta = 
\vert\varphi\vert\circ\theta\circ\varphi\inv
\end{equation}
in \cite{Kaw2} \S2. For each $p \geq 1$, 
{\bf the $p$-th Johnson map} induced by $\theta$
$$
\tau^\theta_p: \Aut(\Fn) \to H^*\otimes H^{\otimes (p+1)}
$$
is defined to be the $p$-th component of the total 
Johnson map $\tau^\theta$, i.e., we have
$$
\tau^\theta(\varphi)\vert_H = 
\sum^\infty_{p=1}\tau^\theta_p(\varphi) \in \Hom(H, \T_2) 
= \prod^\infty_{p=1} H^*\otimes H^{\otimes (p+1)}.
$$
From (1.9) and (2.4) we obtain an integral presentation of
the total Johnson map
\begin{equation}
{\tau^\theta(\varphi)}\inv =
1 +
{\sum}^\infty_{m=1}
\int^{\vert\varphi\vert\circ\theta\circ\varphi\inv}
_{\theta}\overbrace{\eta\eta\cdots\eta}^{\mbox{$m$ times}}
\end{equation}
for any $\varphi \in \Aut(\Fn)$. 
In particular, we have
\begin{equation}
\tau^\theta_1(\varphi) 
= - \int^{\vert\varphi\vert\circ\theta\circ\varphi\inv}
_{\theta}\eta_1.
\end{equation}

Finally we consider a $C^\infty$ path $\theta^t$, $\vert t
\vert\ll 1$, on the space $\Theta_n$. 
Then there exists a $C^\infty$ path $u(t) \in \Lie\,\IA(\T)$ 
such that $u(0) = 0$ and $\theta^t = \exp(u(t))\circ\theta^0$. 
We have $\eta({\thetadot}(0)) =
\stackrel{\centerdot}{u}(0)$ and $\ddt\theta^t(\gamma)
 = \stackrel{\centerdot}{u}(0)\theta^0(\gamma)$ for $\gamma \in
\Fn$. Hence we obtain 
\begin{equation}
\ddt\theta^t(\gamma) =
\eta(\thetadot(0))\theta^0(\gamma)
\end{equation}
for any $\gamma \in \Fn$. Here
$\eta(\thetadot(0)) \in \Lie\,\IA(\T) \cong
H^*\otimes\T_2$ acts on $\theta^0(\gamma) \in \T$ as a
derivation on the algebra $\T$ (2.1).

\section{Stasheff Associahedrons and Twisted Morita-Mumford
Classes}

Let $p \geq 1$ be an integer. The Stasheff associahedron 
$K_{p+1}$ \cite{S} is a finite regular cell complex, 
each of whose cells corresponds to a meaningful way 
of inserting one set of parentheses into the word 
$1 2 3\cdots p+1$ of $p+1$ letters. 
It is homeomorphic to the cube $I^{p-1}$. 
One can observe such a way indicates some twisted 
differential form on the space $\Theta_n$. 
For example, $((123)4)$ indicates $(\eta_2\otimes 1) \eta_1
\in \Omega^2(\Theta_n)\otimes H^*\otimes H^{\otimes 4}$, 
$((1(23))4)$ $(1\otimes \eta_1\otimes 1)(\eta_1\otimes 1)\eta_1
\in \Omega^3(\Theta_n)\otimes H^*\otimes H^{\otimes 4}$, and so on.
Our purpose in this section is to construct a $p$-cochain $Y$ 
in the double cochain complex
$$
Y = Y_p \in C^{*, *}_{(p)} 
:= C^*(K_{p+1}; \left(\Omega^*(\Theta_n)\otimes H^*\otimes
H^{\otimes (p+1)}\right)^{\Aut(\Fn)})
$$
by assembling the Maurer-Cartan forms $\eta_s \in \left(\Omega^*(\Theta_n)\otimes H^*\otimes
H^{\otimes (s+1)}\right)^{\Aut(\Fn)}$, $s \geq 1$.
Here $C^*(K_{p+1})$ means the cellular cochain complex of $K_{p+1}$. 
The Maurer-Cartan formula (2.3) implies $Y_p$ is a $p$-cocycle 
(Proposition 3.1).
Since $K_{p+1}$ is contractible, the cocycle $Y_p$ induces a twisted 
de Rham cohomology class
$$
[Y_p] \in H^p(\left(\Omega^*(\Theta_n)\otimes H^*\otimes
H^{\otimes (p+1)}\right)^{\Aut(\Fn)}).
$$
This is essentially equal to the twisted Morita-Mumford class
$h_p \in H^p(\Aut(\Fn); H^*\otimes H^{\otimes (p+1)})$
\cite{Kaw1}\cite{Kaw2} (Theorem 3.2). \par
In order to construct the $p$-cochain $Y_p$ we fix some notations on 
the Stasheff associahedron $K_{p+1}$. 
Let $S_p$ and $\tilde S_p$ be the sets of all the cells 
and the ordered cells in $K_{p+1}$, respectively. 
We have the canonical projection $\tilde S_p \to S_p$, which is 
clearly a 2-1 map. 
From the definition each element of $S_p$ is regarded as a meaningful 
way of inserting one set of parentheses into the word 
$1 2 \cdots p+1$. For any $\overline{w} \in S_p$ we define the degree
$\deg(\overline{w})$ by the number of pairs of parentheses. 
For example $\deg(((123)4)) = 2$ and $\deg(((1(23))4)) = 3$. 
The dimension of the cell corresponding to $\overline{w}$ is equal to 
$p - \deg(\overline{w})$. 
Following Stasheff \cite{S} we regard $K_{p+1}$ as a subpolyhedron 
of the same dimension in $I^{p-1}$ with the standard orientation. 
Thus we consider $K_{p+1}$ itself as an element of $\tilde S_p$ of degree
$1$.  On the other hand, let $S^0_p \subset S_p$ denote the vertices of 
$K_{p+1}$. We can consider $+w$ and $-w$ as elements of $\tilde S_p$. 
\par
Stasheff introduced the face map $\dd_k(r, s): K_r\times K_s\to 
K_{p+1}$ for $r+s = p+2$ and $1 \leq k \leq r$, to describe the 
boundary $\dd K_{p+1}$. 
The image $\dd_k(r, s)(K_r\times K_s)$ is the cell corresponding 
to $1\cdots k-1(k\cdots k+s-1)k+s\cdots r+s-1$. 
The sign of the map $\dd_k(r, s)$ is $(-1)^{s(r-k)+k}$, namely, 
we have 
$$
\dd K_{p+1} = \sum (-1)^{s(r-k)+k}\dd_k(r, s)(K_r\times K_s)
\in C_{p-2}(K_{p+1})
$$
as cellular chains in $K_{p+1}$. 
The map $\dd_k(r, s)$ should indicate
$\pm 1^{\otimes (k-1)}\otimes \eta_{s-1}\otimes 1^{\otimes (r-k)}
\in \Omega^1(\Theta_n)\otimes H^*\otimes \Hom(H^{\otimes r}, 
H^{\otimes (r+s-1)})$. \par
For simplicity we write
$$
\dd_{a, p', b} := \dd_{a+1}(a+b+1, p'+1): 
K_{a+b+1}\times K_{p'+1} \to K_{a+b+p'+1}
$$
for $a, b \geq 0$ and $p' \geq 1$. 
It indicates $\pm 1^{\otimes a}\otimes \eta_{p'}\otimes 1^{\otimes b}
\in \Omega^1(\Theta_n)\otimes H^*\otimes \Hom(H^{\otimes (a+b+1)}, 
H^{\otimes (a+b+p'+1)})$ and its sign is $(-1)^{(p'+1)b+a+1}$. We have
\begin{equation}
\dd K_{p+1} = \sum_{k, s}(-1)^{(s+1)(p+k)+k+1}\dd_{k, s,
p-k-s}(K_{p-s+1}\times K_{s+1}).
\end{equation}
Looking at an innermost pair of parentheses, we find out 
any oriented cell in $K_{p+1}$ other than $\pm K_{p+1}$ 
is given by $\dd_{a, p', b}(w\times K_{p'+1})$ for some integers 
$a, p', b$ and $w \in \tilde S_{a+b+1}$. 
Clearly $\deg \dd_{a, p', b}(w\times
K_{p'+1}) = \deg w + 1$. In our notation the relations 3(a) and (b) 
in \cite{S} p.278 among the face maps are given by 
\begin{multline}
\dd_{a+k, s, p'-s-k+b}(\dd_{a, p'-s, b}\times 1) = 
\dd_{a, p', b}(1\times \dd_{k, s, p'-k-s}): \\
K_{a+b+1}\times K_{p'-s+1}\times K_{s+1} \to K_{a+b+p'+1}
\end{multline}
and 
\begin{multline}
\dd_{a, p', b+c+p''+1}(\dd_{a+b+1, p'', c}\times 1)(1\times T) = 
\dd_{a+b+p'+1, p'', c}(\dd_{a, p', b+c+1}\times 1): \\
K_{a+b+c+1}\times K_{p'+1}\times K_{p''+1} \to K_{a+b+c+p'+p''+1},
\end{multline}
respectively.
Here $T: K_{p'+1}\times K_{p''+1} \to K_{p''+1}\times K_{p'+1}$ 
is the switch map.\par
Now we can define the map
$$
Y: {\coprod}_{p\geq 1}\tilde S_p \to {\bigoplus}_{p\geq 1}
\left(\Omega^*(\Theta_n)\otimes H^*\otimes
H^{\otimes (p+1)}\right)^{\Aut(\Fn)}
$$
inductively on $\deg w$ by 
\begin{eqnarray}
&& Y(\pm K_{p+1}) := \pm \eta_p, \quad \text{and}\\
&& Y(\dd_{a, p', b}(w\times K_{p'+1})) := (-1)^{a + (p'+1)b + \deg w}
(1^{\otimes a}\otimes \eta_{p'}\otimes 1^{\otimes b})Y(w). \nonumber
\end{eqnarray}
An innermost pair of parentheses is not necessarily unique. 
We have to prove the definition (3.4) is independent of the 
choice of innermost pairs.\par
It suffices to compute the sign related to the twisted $2$-form 
$1^{\otimes a}\otimes \eta_{p'}\otimes 1^{\otimes b}\otimes \eta_{p''}\otimes
1^{\otimes c}$ in two different ways. Let $w$ be an element 
of $\tilde S_{a+b+c+1}$. By a straightforward computation 
one can obtain
\begin{eqnarray*}
&& (-1)^{(p'-1)(p''-1)}Y(\dd_{a, p', b+c+p''+1}(\dd_{a+b+1, p'', c}\times
1)(w\times K_{p''+1}\times K_{p'+1}))\\
&=& (-1)^{(p'+1)(b+c)+(p''+1)c +b}(1^{\otimes a}\otimes \eta_{p'}\otimes
1^{\otimes b}\otimes \eta_{p''}\otimes 1^{\otimes c})Y(w)\\
&=& Y(\dd_{a+b+p'+1, p'', c}(\dd_{a, p', b+c+1}\times 1)(w\times
K_{p'+1}\times K_{p''+1})).
\end{eqnarray*}

The restriction $Y_p := Y\vert_{\tilde S_p}: \tilde S_p \to 
\left(\Omega^*(\Theta_n)\otimes H^*\otimes
H^{\otimes (p+1)}\right)^{\Aut(\Fn)}$ can be regarded as a $p$-cochain 
of the double complex 
$$
C^{*, *}_{(p)} 
= C^*(K_{p+1}; \left(\Omega^*(\Theta_n)\otimes H^*\otimes
H^{\otimes (p+1)}\right)^{\Aut(\Fn)}).
$$
The Maurer-Cartan formula (2.3) implies

\begin{proposition}
$$
(dY)(w) = (\delta Y)(w)
$$
for any $p \geq 1$ and $w \in \tilde S_p$. 
In other words, $Y_p \in C^{*, *}_{(p)}$ is a $p$-cocycle.
\end{proposition}
\begin{proof} We prove it by induction on $\deg w$. 
By (2.3) and (3.1) we have 
\begin{eqnarray*}
&& (dY)(K_{p+1}) = d\eta_p = \sum_{k, s}(1^{\otimes k}\otimes
\eta_{s}\otimes 1^{\otimes (p-s-k)})\eta_{p-s}\\
&=& \sum_{k, s}(-1)^{(s+1)(p+k)+k+1}Y(\dd_{k, s,
p-k-s}(K_{p-s+1}\times K_{s+1}))\\
&=& Y(\dd K_{p+1}) = (\delta Y)(K_{p+1}).
\end{eqnarray*}
Assume $w \in \tilde S_{a+b+1}$ satisfies the proposition. 
We have 
\begin{eqnarray*}
&& dY(\dd_{a, p', b}(w\times K_{p'+1})) \\
&=& (-1)^{a + (p'+1)b + \deg w}
d((1^{\otimes a}\otimes \eta_{p'}\otimes 1^{\otimes b})Y(w))\\
&=& (-1)^{a + (p'+1)b + \deg w}
(1^{\otimes a}\otimes d\eta_{p'}\otimes 1^{\otimes b})Y(w)\\
&& + (-1)^{a + (p'+1)b + \deg w+1}
(1^{\otimes a}\otimes \eta_{p'}\otimes 1^{\otimes b})dY(w).\\
\end{eqnarray*}
Then from the Maurer-Cartan formula (2.3), (3.2) and (3.1) 
follows
\begin{eqnarray*}
&&(-1)^{p'b}(1^{\otimes a}\otimes d\eta_{p'}\otimes 1^{\otimes b})Y(w)\\
&=& (-1)^{p'b}\sum_{k, s}(1^{\otimes (a+k)}\otimes \eta_{s}\otimes
1^{\otimes (p'-s-k+b)})(1^{\otimes a}\otimes \eta_{p'-s}\otimes
1^{\otimes b})Y(w)\\
&=& \sum_{k, s}(-1)^{k+(s+1)(p'-k)+1}
Y(\dd_{a+k, s, p'-s-k+b}(\dd_{a, p'-s, b}\times 1)(w\times K_{p'-s+1}\times
K_{s+1}))\\
&=& \sum_{k, s}(-1)^{k+(s+1)(p'-k)+1}
Y(\dd_{a, p', b}(1\times \dd_{k, s, p'-k-s})(w\times K_{p'-s+1}\times
K_{s+1}))\\
&=& Y(\dd_{a, p', b}(w\times \dd K_{p'+1})).
\end{eqnarray*}
By the inductive assumption we have
\begin{eqnarray*}
&& (1^{\otimes a}\otimes \eta_{p'}\otimes 1^{\otimes b})dY(w)
= (1^{\otimes a}\otimes \eta_{p'}\otimes 1^{\otimes b})(\delta Y)(w)\\
&=& (1^{\otimes a}\otimes \eta_{p'}\otimes 1^{\otimes b}) Y(\dd w)
= (-1)^{a + (p'+1)b + \deg w +1}Y(\dd_{a, p', b}((\dd w)\times K_{p'+1})).
\end{eqnarray*}
Therefore we obtain 
\begin{eqnarray*}
&& dY(\dd_{a, p', b}(w\times K_{p'+1})) \\
&=& (-1)^{a+b-\deg w}Y(\dd_{a, p', b}(w\times \dd K_{p'+1})
+ Y(\dd_{a, p', b}((\dd w)\times K_{p'+1}))\\
&=& Y(\dd_{a, p', b}\dd(w\times K_{p'+1})) 
= Y(\dd(\dd_{a, p', b}(w\times K_{p'+1})))\\
&=& (\delta Y)(\dd_{a, p', b}(w\times K_{p'+1})).
\end{eqnarray*}
This completes the induction.
\end{proof}

We define 
$$
w^0_p := ((((12)3)\cdots)p+1) \in S^0_p \subset \tilde S_p.
$$
Then $\deg w^0_p = p$ and $w^0_p = \dd_{0, 1, p-1}(w^0_{p-1}\times K_2)$. 
Hence we have 
\begin{eqnarray}
Y_p(w^0_p) &=& (-1)^{\frac12p(p-1)}(\eta_1\otimes 1^{\otimes (p-1)})
\cdots(\eta_1\otimes 1)\eta_1\nonumber\\
&=& (-1)^{\frac12p(p-1)}{\bb_p}_*({\eta_1}^{\otimes p}) \in
\Omega^p(\Theta_n)\otimes H^*\otimes H^{\otimes (p+1)},
\end{eqnarray}
where $\bb_p: (H^*\otimes H^{\otimes 2})^{\otimes p} 
\to H^*\otimes H^{\otimes(p+1)}$ is a $\GL(H)$-homomorphism 
indicated by $w^0_p$ \cite{Kaw2} \S4. In ibid.\! Theorem 4.1 
we proved 
\begin{equation}
h_p = {\bb_p}_*([\tau^\theta_1]^{\otimes p}) 
\in H^p(\An; H^*\otimes H^{\otimes(p+1)})
\end{equation}
for any $p \geq 1$. 
Here $\theta$ is any Magnus expansion and $h_p$ is a cohomology class 
introduced in ibid.\! \S4. 
The restriction to the mapping class group $\Mgone$ of genus $g$ 
with $1$ boundary component is essentially equal to the 
$(0, p+2)$-twisted Morita-Mumford class \cite{Kaw1}
$$
(p+2)!\,h_p\vert_{\Mgone} = m_{0, p+2} 
\in H^p(\Mgone; H^*\otimes H^{\otimes (p+1)})
$$
\cite{Kaw2} (5.8). 
As was shown in \cite{Mo3}, \cite{KM1} and \cite{KM2}, each of the
Morita-Mumford classes
$e_i$ is obtained  by contracting the coefficients of the twisted ones 
in a suitable way using the intersection form on the homology group 
of the surface. \par
Let $M$ be a connected $C^\infty$ manifold, 
$\widehat{f}: \pi_1(M) \to
\Aut(\Fn)$ a group homomorphism, 
and $f: \widetilde{M} \to \Theta_n$ an
equivariant $C^\infty$ map 
of the universal covering space $\widetilde{M}$ of $M$ 
into $\Theta_n$ with
respect to $\widehat{f}$.
Then we obtain the pullback
$$
f^*Y_p \in C^*(K_{p+1}; \Omega^*(M; H^*\otimes H^{\otimes (p+1)}))
$$
for any $p \geq 1$. 
Hence we regard $H^*\otimes H^{\otimes (p+1)}$ as a flat vector 
bundle over $M$ through the homomorphism $\widehat{f}$. 
Since $K_{p+1}$ is contractible, an isomorphism 
\begin{equation}
H^*(C^*(K_{p+1}; \Omega^*(M; H^*\otimes H^{\otimes (p+1)}))) 
= H^*(M;  H^*\otimes H^{\otimes (p+1)})
\end{equation}
holds and maps $[f^*Y_p]$ to the twisted de Rham class 
$(-1)^{\frac12p(p-1)}{\bb_p}_*({[\eta_1]}^{\otimes p})$. 
From (2.5) and (3.6) we have $h_1 = -[\eta_1] \in 
H^1(M; H^*\otimes H^{\otimes 2})$. Therefore we obtain
\begin{theorem}
$$
f^*[Y_p] = (-1)^{\frac12p(p+1)}{\widehat{f}}^*h_p \in 
H^p(M; H^*\otimes H^{\otimes(p+1)})
$$
under the isomorphism (3.7).
\end{theorem}
Let $\moduli$ be the moduli space of 
triples $(C, P_0, v)$ of genus $g \geq 1$, 
where $C$ is a compact Riemann surface 
of genus $g$, $P_0 \in C$, and $v$ a non-zero tangent 
vector of $C$ at $P_0$. It is known that $\moduli = B\Mgone$. 
In Part II we will construct an equivariant real analytic map 
$\theta: \widetilde{\moduli} \to \Theta_{2g}$ we call {\bf
the harmonic Magnus expansion} (\S\S4-6) and 
compute the pullback $\theta^*\eta_p$ in an 
explicit way (\S\S7-8).

\part{Harmonic Magnus Expansions}
\section{Fundamental Group with a Tangential
Basepoint}

Let $g \geq 1$ be a positive integer. 
In \S\S4-6 we construct {\bf the harmonic Magnus 
expansion} of a {\bf triple} $(C, P_0, v)$ {\bf of genus} $g$
$$
\theta^{(C, P_0, v)}: \pi_1(C, P_0, v) \to 
1+\T_1(H_1(C; \bR)),
$$
where $C$ is a compact Riemann surface of genus $g$, 
$P_0 \in C$ and $v \in T_{P_0}C\setminus\{0\}$. 
$\pi_1(C, P_0, v)$ is the fundamental group of 
the triple $(C, P_0, v)$, 
which will be defined in this section as the fundamental 
group of the complement $C\setminus \{P_0\}$ with the 
tangential basepoint $v$. 
Moreover $\T = \T(H_1(C; \bR))$ is the completed 
tensor algebra generated by the first real homology 
group $H_1(C; \bR)$ of the Riemann surface $C$
$$
\T(H_1(C; \bR)) := 
{\prod}^\infty_{m=0}H_1(C; \bR)^{\otimes m}.
$$
The two-sided ideal $\T_p(H_1(C; \bR))$ for each $p \geq 1$ 
is defined as in \S2. \par
Up to \S6 we study a {\bf single} triple $(C, P_0, v)$. 
So we write simply $H = H_1(C; \bR)$, $\T = \T(H_1(C; \bR))$ 
and $\T_p = \T_p(H_1(C; \bR))$. \par
Now we define the fundamental group $\pi_1(C, P_0, v)$ of the 
triple $(C, P_0, v)$. 
We denote by $\frak L(C, P_0, v)$ the set of all the piecewise 
$C^\infty$ map $\ell: [0, 1] \to C$ satisfying the conditions 
\begin{eqnarray}
&& \ell(]0, 1[) \subset C \setminus \{P_0\},\\
&& \ell(0) = \ell(1) = P_0, \quad\text{and}\nonumber\\
&& \frac{d\ell}{dt}(0) = -\frac{d\ell}{dt}(1) = v.\nonumber 
\end{eqnarray}
For $\ell_0$ and $\ell_1 \in \frak L(C, P_0, v)$ 
we define $\ell_0 \sim \ell_1$ if there exists a piecewise 
$C^\infty$ map $L: [0, 1] \times [0, 1] \to C$ such that 
\begin{eqnarray}
&& L(]0, 1[\times [0, 1]) \subset C \setminus \{P_0\},\\
&& L(0, s) = L(1, s) = P_0, \quad \forall s \in [0, 1],\nonumber\\
&& \frac{\partial L}{\partial t}(0, s) = -\frac{\partial
L}{\partial t}(1, s) = v,
\quad \forall s \in [0, 1], \quad\text{and}\nonumber\\
&& L(t, 0) = \ell_0(t), \quad L(t, 1) = \ell_1(t),
\quad \forall t \in [0, 1].\nonumber
\end{eqnarray}
Clearly $\sim$ is an equivalent relation. 
We call the quotient set 
$$
\pi_1(C, P_0, v) := \frak L(C, P_0, v)/\sim
$$
{\bf the fundamental group of the triple} 
$(C, P_0, v)$ (with the tangential basepoint $v$).
The set $\pi_1(C, P_0, v)$ has a natural group structure 
isomorphic to the free group of rank $2g$, $F_{2g}$. 
Its abelianization is naturally isomorphic to the 
first integral homology group of $C\setminus\{P_0\}$ 
and that of $C$
$$
\pi_1(C, P_0, v)\abel =  H_1(C \setminus \{P_0\};\, \bZ) = H_1(C;\, \bZ).
$$
\par

We choose a symplectic generator 
$\{x_i\}^{2g}_{i=1}$ of the fundamental group 
$\pi_1(C, P_0, v)$. 
It gives an isomorphism $F_{2g} \overset\cong\to
 \pi_1(C, P_0, v)$. A negative loop arround $P_0$ gives 
a word
$$
w_0 := {\prod}^g_{i=1} x_ix_{g+i}{x_i}\inv{x_{g+i}}\inv.
$$
The first real homology group has the intersection number 
$\cdot$ satisfying
$$
X_i\cdot X_{g+j} = \delta_{i, j}, \quad\text{and}\quad
X_i\cdot X_j = X_{g+i}\cdot X_{g+j} = 0
$$
for $1 \leq i, j \leq g$, where $X_i = [x_i] \in H$ 
as in \S2. The intersection form $I$ is equal to  
$$
I = {\sum}^g_{i=1}(X_iX_{g+i} -  X_{g+i} X_i)
\in H^{\otimes 2}. 
$$

The Poincar\'e duality $\vartheta = \cap [C] : H^* = H^1(C; \bR)
\to H = H_1(C; \bR)$ and its inverse $\vartheta\inv$ are 
given by 
\begin{eqnarray}
&& \vartheta: H^* \overset\cong\to H, \quad
\xi \mapsto -(\xi\otimes 1_H)(I), \quad \text{and}\\
&& \vartheta\inv: H \overset\cong\to H^*, \quad
Z \mapsto Z\cdot,\nonumber
\end{eqnarray}
respectively. 
For details, see \cite{Kaw2} \S5. Throughout Part II we identify 
$H$ and its dual $H^*$ by the isomorphism $\vartheta$, which is
equivariant under the action of the mapping class group.\par
Any element $u = \sum^{2g}_{i=1} u_iX_i \in \T\otimes H = \T_1$,
$u_i \in \T$, induces a derivation $\inter(u)$ on the algebra
$\T$. In fact, for any $Z_j \in H$, $1 \leq j \leq m$, we define
$$
\inter(u)\left(Z_1Z_2\cdots Z_m\right) :=
\inter(u)(Z_1)Z_2\cdots Z_m 
+\cdots + Z_1\cdots Z_{m-1}\inter(u)(Z_m)
$$
and $\inter(u)(Z_j) := {\sum}^{2g}_i u_i(X_i\cdot Z)$. 
Thus we obtain a linear map 
$$
\inter: \T_1 \to \operatorname{Der}(\T).
$$
The $m$-th symmetric group $\frak S_m$ acts on 
$H^{\otimes m}$ by permuting the components. 
We define a linear automorphism $\varepsilon$ of $\T_1$ by 
\begin{equation}
\varepsilon\vert_{H^{\otimes m}} := \begin{pmatrix}
1& 2& \cdots & m-1 & m\\
2& 3& \cdots & m & 1
\end{pmatrix}.
\end{equation}
It is easy to show 
\begin{equation}
\inter(u)(I) = \varepsilon u - u,
\end{equation}
which seems to be related to \cite{Mo4} Proposition 4.6, p.366. \par
As in (2.2) the space $\Hom(H, \T_2) = H^*\otimes \T_2$ can be 
regarded as a space of derivations on the algebra $\T$. 
From (4.3) this action coincides with the composite
\begin{equation}
\Hom(H, \T_2) = H^*\otimes \T_2
\overset{\vartheta\otimes 1}\longrightarrow
H\otimes \T_2 \overset{\varepsilon\inv}\longrightarrow
\T_3 \overset{\inter}\longrightarrow
\operatorname{Der}(\T). 
\end{equation}

\section{The Connection $1$-form $\omega$}

In \S 6 we will construct the harmonic Magnus expansion 
of any triple $(C, P_0, v)$ of genus $g \geq 1$ 
by using Chen's iterated integrals \cite{C}, 
whose integrand is a real-valued connection $1$-form
$$
\omega \in A^1_{\bR}(C) \,\wotimes\,\T_1 
= A^1_{\bR}(C) \,\wotimes\,\T_1(H_1(C; \bR))
= {\prod}^\infty_{m=1}A^1_{\bR}(C)\otimes H_1(C; \bR)^{\otimes m}
$$
singular at the point $P_0$, 
constructed in a canonical way in this section. 
Here we denote by $A^1_{\bR}(C)$ the real-valued
$1$-currents on $C$. \par

We begin by recalling some results on a Green operator 
for the compact Riemann surface $C$ of genus $g$. 
For $q \geq 0$ we denote by $A^q(C)$ the complex-valued
$q$-currents on $C$. 
The Hodge $*$-operator $*: \left(T^*_\bR C\right)\otimes\bC 
\to \left(T^*_\bR C\right)\otimes \bC$ on the cotangent 
bundle of $C$ depends only on the complex structure of $C$.
The $-\sqrt{-1}$-eigenspace is the holomorphic cotangent bundle
$T^*C$, and the $\sqrt{-1}$-eigenspace the antiholomorphic
cotangent bundle $\overline{T^*C}$. 
We have an exact sequence
\begin{equation}
0 \to \bC \to A^0(C) \overset{d*d}\to A^2(C)
\overset{\int_C}\to \bC \to 0
\label{5.1}\end{equation}
The vector space $\bC$ on the left side means the constant
functions. A {\bf Green operator}
$$
\Phi = \Phi^{(C, P_0)}: A^2(C) \to A^0(C)/\bC
$$
is defined by
\begin{equation}
d* d\Phi \Omega = \Omega - \left(\int_C\Omega\right)\delta_0
\label{5.2}
\end{equation}
for any $\Omega \in A^2(C)$.
Here $\delta_0 = \delta_{P_0}: C^\infty(C) \to \bC$,
$f \mapsto f(P_0)$, is the delta current on $C$ at 
the point $P_0$.
\par
The operator $\Phi$ gives the Hodge
decomposition of the $1$-currents 
\begin{eqnarray}
\varphi = \harmonic\varphi + d\Phi d\ast \varphi + \ast d\Phi d \varphi
\label{5.3}
\end{eqnarray}
for any $\varphi \in A^1(C)$, where $\harmonic: A^1(C) \to A^1(C)$ is the
harmonic projection on the $1$-currents on $C$.\par

The Hodge $*$-operator decomposes the space $A^1(C)$ into 
the $\pm \sqrt{-1}$-eigenspaces
$$
A^1(C) = A^{1, 0}(C) \oplus A^{0, 1}(C),
$$ 
where $A^{1, 0}(C)$ is the $-\sqrt{-1}$-eigenspace and 
$A^{0, 1}(C)$ the $\sqrt{-1}$-eigenspace. 
Throughout Part II we denote by $\varphi'$ and $\varphi''$ 
the $(1, 0)$- and the $(0, 1)$-parts of $\varphi \in A^1(C)$, 
respectively, i.e., 
$$
\varphi = \varphi' + \varphi'', \quad 
*\varphi = -\sqrt{-1}\varphi' + \sqrt{-1}\varphi''.
$$
If $\varphi$ is harmonic, then $\varphi'$ is holomorphic and 
$\varphi''$ anti-holomorphic.\par

Our purpose is to define a connection $1$-form
$$
\omega \in A^1(C) \,\wotimes\,\T_1
$$
for any pointed compact Riemann surface $(C, P_0)$
in a canonical way.
If $\varphi \in A^1(C)\wotimes\T = \prod^\infty_{m=0}A^1(C)\otimes
H^{\otimes m}$ is a
$\T$-valued
$q$-current,  then we denote by $\varphi_{(m)}$ the $m$-th homogeneous term
of $\varphi$
$$
\varphi = {\sum}^\infty_{m=0}\varphi_{(m)}, \quad \varphi_{(m)} \in
A^1(C)\otimes H^{\otimes m}.
$$
\par
In \S7 we have to evaluate the function $\Phi\Omega$ at the point
$P_0$. Of course, $\Phi\Omega$ has not necessarily its own value
at the singular point $P_0$ for any $2$-current $\Omega$. 
Hence we need to introduce some appropriate function spaces.
Our consideration is based on the following result
which Ahlfors and Bers used to prove
Riemann's mapping theorem for variable metrics
(\cite{AB} Lemmas 1 and 3, p.386).
\begin{theorem}\label{Theorem 5.1} Let $p > 2$.
If $f$ is a distribution defined on an open subset in $\bC$,
and
its derivative $f_\zbar$ is locally $L^p$, then $f$ itself
is a H\"older continuous function of class
$C^{0+(1-2/p)}$.
\end{theorem}
We define the spaces $E^1(C)$, $E^0(C)$ and $E^0_0(C)$ by
\begin{eqnarray}
&& E^1(C) := \bigcap_{2 < p < \infty}L^p\left(C; \, (T^*_\bR
C)\otimes\bC\right), \label{5.4}\\
&& E^0(C) := \bigcap_{2 < p < \infty}C^{0+(1-2/p)}
\left(C; \, \bC\right), \quad \text{and}
\nonumber\\
&& E^0_0(C, P_0) := \{f \in E^0(C); \, f(P_0) = 0\},
\nonumber
\end{eqnarray} 
respectively. Then, for any $\varphi \in
E^1(C)$, $\Phi d\varphi$ is contained in $E^0(C)$ up to
translation of constant functions. In fact, the 
$(0, 1)$-part satisfies 
$$
2\sqrt{-1}\overline{\partial}\Phi d \varphi''
= (d\Phi d* + *d\Phi d) \varphi''
= (1-\harmonic) \varphi'' \in E^1(C),
$$
which implies $\Phi d \varphi'' \in E^0(C)$ by Theorem 5.1.
Similarly we have $\Phi d \varphi' \in E^0(C)$. \par
Hence we may regard $\Phi d$ as an operator
$E^1(C) \to E^0(C)/\bC$, and may normalize it by
\begin{equation}
(\Phi d\varphi)(P_0) = 0
\label{5.5}\end{equation}
for any $\varphi \in E^1(C)$, namely, we consider it a linear map
$\Phi d: E^1(C) \to E^0_0(C, P_0)$.\par

Now we define the $m$-th homogenous term
$\omega_{(m)} \in A^1(C)\otimes H^{\otimes m}$ of
the connection form $\omega$ by induction on $m$.
The first term $\omegaone \in A^1(C)\otimes H$ comes 
from the real harmonic $1$-forms on $C$. 
As in \S4 we choose a symplectic basis $\{X_i\}^{2g}_{i=1}$ 
of $H$. Let $\{\xi_i\}^{2g}_{i=1}$ be the basis of the real 
harmonic $1$-forms on $C$ whose cohomology classes $[\xi_i]$, 
$1 \leq i \leq 2g$, form the dual basis of $\{X_i\}^{2g}_{i=1}$, i.e.,
$$
\int_{X_j}\xi_i = \delta_{i, j}
$$
for $1 \leq i, j \leq 2g$.
We have $\vartheta[\xi_i] = -X_{g+i}$, $\vartheta[\xi_{g+i}] =
X_i$, 
\begin{equation}
\int_C\xi_i\wedge\xi_{g+j} = \delta_{i, j}, 
\quad\text{and}\quad
\int_C\xi_i\wedge\xi_j = \int_C\xi_{g+i}\wedge\xi_{g+j} = 0
\label{5.6}\end{equation}
for $1 \leq i, j \leq g$. We define 
\begin{equation}
\omegaone := \sum^g_{i=1} \xi_iX_i + \xi_{g+i}X_{g+i}. 
\label{5.7}\end{equation}
Clearly it is independent of the choice of a basis 
$\{X_i\}^{2g}_{i=1}$, and satisfies 
\begin{equation}
\int_\ell\omegaone = [\ell] \in H = H_1(C; \bR) 
\label{5.8}
\end{equation}
for any $\ell \in \frak L(C, P_0, v)$. 
The harmonic projection $\harmonic: A^1(C)\wotimes\T \to
A^1(C)\wotimes\T$ is given by
\begin{equation}
\harmonic\varphi
= - \left(\int_C\varphi\wedge\omegaone\right)\cdot\omegaone
\label{5.9}
\end{equation}
for any $\varphi \in A^1(C)\wotimes\T$. In fact, if $\harmonic
\varphi = \sum^g_{i=1} x_i\xi_i + y_i\xi_{g+i}$, $x_i, y_i \in
\T$, then the right-hand side of (5.9) is 
\begin{eqnarray*}
&& -\left(\int_C \varphi \wedge\omegaone\right)\cdot\omegaone
= -\left(\int_C (\harmonic\varphi) \wedge\omegaone\right)\cdot\omegaone
\\
&=& -\left(\int_C \left(\sum^g_{i=1} x_i\xi_i +
y_i\xi_{g+i}\right)
\wedge\left(\sum^g_{j=1}\xi_jX_j +
\xi_{g+j}X_{g+j}\right)\right)\cdot\omegaone
\\
&=& -\left(\sum^g_{i=1} x_iX_{g+i} - y_iX_i\right)\cdot
\left(\sum^g_{j=1}\xi_jX_j + \xi_{g+j}X_{g+j}\right)
=  \sum^g_{i=1} x_i\xi_i + y_i\xi_{g+i}
= \harmonic \varphi.
\end{eqnarray*}
Moreover we have 
\begin{equation}
\int_C\omegaone\wedge\omegaone = {\sum}^g_{i=1}X_iX_{g+i} - X_{g+i}X_i = I. 
\label{5.10}\end{equation}
Now we define $\omega_{(m)}$ for $m \geq 2$ inductively by
\begin{equation}
\omega_{(m)} := \ast d\Phi(\omega\wedge\omega)_{(m)}
=\ast d\Phi\left({\sum}^{m-1}_{p=1}\omega_{(p)}
\wedge\omega_{(m-p)}\right).
\label{5.11}
\end{equation}
Then immediately we have
\begin{eqnarray}
&& d\omega = \omega\wedge\omega -
I\,\delta_0,\label{5.12}\\
&& * d\Phi(\omega\wedge\omega) = \omega - \omegaone.
\label{5.13}
\end{eqnarray}
\par
Next we study the regularity of the connection form $\omega_{(m)}$. 
We use
the following classical result, whose proof one 
can find in Bers'
lecture note \cite{B}.
\begin{theorem}\label{Theorem 5.2}{\rm (\cite{L}\cite{Ko})} Let $f$ be a
distribution defined on an open subset in $\bC$. Suppose $f_\zbar$ is
H\"older continuous of class $C^{l+\alpha}$, 
$\alpha \in ]0, 1[$, $l \in \bZ_{\geq 0}$. 
Then $f$ is of class $C^{l+1+\alpha}$.
\end{theorem}

Clearly the harmonic form $\omegaone$ is of class $C^\infty$
on the whole $C$. Theorem 5.2 implies $\omega_{(m)}$
is of class $C^\infty$ on $C - \{P_0\}$ for any $m \geq 2$.
Hence it suffices to study the regularity near the
point $P_0$.\par

Immediately from Theorems 5.1 and 5.2 we have

\begin{corollary}\label{Corollary 5.3} Let $f$ be a distribution
defined on an open subset in $\bC$. Suppose $d*df (= 
2\sqrt{-1}\dd\dc f)$ is locally $L^p$ for any 
$p \in ]2, +\infty[$. Then $f$ is H\"older continuous 
of class $C^{1+\alpha}$ for any $\alpha \in ]0, 1[$, 
and $*df$ of class $C^{0+\alpha}$.
\end{corollary}

We fix a complex coodinate $z$ of $C$ centered at $P_0$. 
In the polar coordinates
$$
z = re^{\sqrt{-1}\theta}, \quad r \in \bR_{\geq 0}, \, \, 
\theta = \arg z \in \bR/2\pi\bZ,
$$
the $1$-current 
\begin{equation}
\frac1{2\pi}d\theta 
= \frac1{2\pi}d\arg z
= \frac1{4\pi\sqrt{-1}}\left(\frac{dz}z - \frac{d\zbar}\zbar\right) 
= *d\left(\frac1{2\pi}\log\vert z\vert\right)\label{5.14}
\end{equation}
satisfies
\begin{equation}
\delta_0 = d\left(\frac1{2\pi}d\theta\right) =
d\left(\frac1{2\pi}d\arg z\right) \label{5.15}
\end{equation}
near $P_0$. In fact, we have 
\begin{eqnarray*}
&& 
\int_{\vert z\vert \ll1}d\left(\frac1{2\pi}d\theta\right)f = 
\int_{\vert z\vert \ll1}\left(\frac1{2\pi}d\theta\right)\wedge df\\
&=& \lim_{\epsilon\downarrow0}\int_{\epsilon\leq \vert z\vert
\ll1}\left(\frac1{2\pi}d\theta\right)\wedge df =
\lim_{\epsilon\downarrow0}\int_{\vert z\vert
=\epsilon}\frac1{2\pi}f(z)d\theta = f(P_0)
\end{eqnarray*}
for any $C^\infty$ function $f$ with a compact support near $P_0$. 
From (5.12), (5.13) and (5.15) the $1$-current
$$
\widetilde{\omega}_{(2)} := \omega_{(2)} - \frac1{2\pi}Id\theta
$$
is of class $C^\infty$ near $P_0$.\par
For the rest of this section we prove
\begin{proposition}\label{Proposition 5.4} 
Let $z = re^{\sqrt{-1}\theta}$ be a complex coordinate of $C$ centered at 
$P_0$ as above. Then 
\begin{itemize}
\item[(1)] $\omegaone$ is of class $C^\infty$ near $P_0$.
\item[(2)] $\widetilde{\omega}_{(2)}$ is of class $C^\infty$ near $P_0$.
\item[(3)] For $m \geq 3$ we have some polynomial $u_m(s) \in H^{\otimes
n}\otimes \bC[s]$ of degree $\leq \frac12(m-1)$ such that 
$$
\omega_{(m)} - *d\left(u_m(\log\vert z\vert)z+ \overline{u_m(\log\vert
z\vert)}\zbar\right)
$$
is H\"older continuous of class $C^{0+\alpha}$ near $P_0$ for any $\alpha \in
]0, 1[$.
\item[(4)] For $m \geq 3$ we have some constant $C_m > 0$ such that 
$$
\left\vert\omega_{(m)}\left(\frac{d}{dz}\right)\right\vert = 
\left\vert\omega_{(m)}\left(\frac{d}{d\zbar}\right)\right\vert \leq
C_m \left\vert \log\vert z\vert\right\vert^{[(m-1)/2]}
$$
near $P_0$. 
\end{itemize}
\end{proposition}

\begin{proof} (1) is clear. (2) was already shown. To prove (3) and (4), 
we need the following, which is obtained by straightforward computation.
\begin{lemma}\label{Lemma 5.5} 
For any polynomial $c(s) \in \bC[s]$ we have 
$$
*d(c(\log\vert z\vert)z) 
= -\sqrt{-1}\left(c(\log\vert z\vert) + \frac12c'(\log\vert z\vert)\right) dz
+ \sqrt{-1}c'(\log\vert z\vert)\frac{z}{2\zbar}d\zbar,
\leqno{(1)} 
$$
$$
d*d(c(\log\vert z\vert)z) 
= \sqrt{-1}\left(c'(\log\vert z\vert) 
+ \frac12c''(\log\vert z\vert)\right) dz\wedge 
\left(\frac1\zbar d\zbar\right),
\leqno{(2)} 
$$
$$
\left(\frac1{2\pi}d\theta\right)\wedge*d(c(\log\vert z\vert)z) 
= \frac{-1}{4\pi}\, c(\log\vert z\vert) dz\wedge 
\left(\frac{1}{\zbar}d\zbar\right).
\leqno{(3)}
$$
\end{lemma}
Now we prove (3) and (4) by induction on $m \geq 3$. 
Our proof for $\omega_{(3)}$ and $\omega_{(4)}$ is 
slightly different from that for $\omega_{(m)}$ for $m \geq 5$. 
Here we should remark 
\begin{equation}
\int_C\omega_{(q)}\wedge\omega_{(m-q)} = 0
\label{5.16}\end{equation}
for any $q$, $1 \leq q \leq m-1$.\par
First we study $\omega_{(3)}$. The value of the real form 
$\omegaone$ at $P_0$, $\omegaone(P_0)$, is given by 
$$
\omegaone(P_0) = w_1(dz)_{P_0} + \overline{w_1}(d\zbar)_{P_0}
$$
for some $w_1 \in H\otimes\bC$. 
From Lemma 5.5(2) we have 
$$
w_1dz\wedge \left(\frac1{2\pi}d\theta\right) 
= \frac{\sqrt{-1}}{4\pi} w_1dz\wedge \frac{d\zbar}\zbar
= \frac{1}{4\pi}d*d(w_1z\log\vert z\vert).
$$
Hence, if we define 
$$
u_3(s) := \frac{s}{4\pi}(w_1I-Iw_1) \in H^{\otimes 3}\otimes \bC[s],
$$
then 
$$
\omegaone\wedge\left(\frac1{2\pi}I d\theta\right) +
\left(\frac1{2\pi}Id\theta\right)\wedge\omegaone 
- d*d\left(u_3(\log \vert z \vert)z + \overline{u_3(\log \vert z
\vert)}\zbar\right)
$$
and 
\begin{eqnarray*}
&& \omegaone\wedge\omega_{(2)} + \omega_{(2)} \wedge\omegaone 
-  d*d\left(u_3(\log \vert z \vert)z + \overline{u_3(\log \vert z
\vert)}\zbar\right)
\\
&=&d*d\left(\Phi(\omega\wedge\omega)_{(3)} 
- \left(u_3(\log \vert z \vert)z + \overline{u_3(\log \vert z
\vert)}\zbar\right)\right)
\end{eqnarray*}
are bounded near $P_0$. From Corollary 5.3
$$
\omega_{(3)} -*d\left(u_3(\log \vert z \vert)z + 
\overline{u_3(\log \vert z\vert)}\zbar\right)
$$
is H\"older continuous of class $C^{0+\alpha}$ for any $\alpha \in ]0, 1[$.
Clearly $\deg u_3(s) \leq 1$. 
From Lemma 5.5(1) we have some constant $\widehat{C_3} > 0$ such that 
$$
\left\vert *d(u_3(\log \vert z \vert)z) \left(\frac{d}{dz}\right)
\right\vert \leq \widehat{C_3}\vert\log\vert z\vert\vert
$$
near $P_0$. This proves (3) and (4) for the case $m = 3$. \par
Next we consider the case $m = 4$. The $2$-current 
$$
\omega_{(2)}\wedge\omega_{(2)} - 
\widetilde{\omega}_{(2)}\wedge \frac1{2\pi}Id\theta
-\frac1{2\pi}Id\theta\wedge \widetilde{\omega}_{(2)}
$$
is of class $C^\infty$ near $P_0$. 
The value of $\widetilde{\omega}_{(2)}$ at $P_0$ is given by 
$$
\widetilde{\omega}_{(2)}(P_0) = w_2\left(dz \right)_{P_0} +
\overline{w_2}\left(d\zbar \right)_{P_0} 
$$
for some $w_2 \in H^{\otimes 2}\otimes \bC$. 
By a similar argument to that for $\omega_{(3)}$, if we define 
$$
u_4(s) := \frac{s}{4\pi}\left(w_2 I -Iw_2 \right)
\in H^{\otimes 4}\otimes\bC[s],
$$
then 
$$
\omega_{(2)}\wedge\omega_{(2)} - d*d\left(u_4(\log \vert z \vert)z +
\overline{u_4(\log \vert z\vert)}\zbar\right)
$$
is bounded near $P_0$. On the other hand, 
by the assertion (4) for $\omega_{(3)}$, the $2$-current 
$\omegaone\wedge\omega_{(3)} + \omega_{(3)}\wedge\omegaone$ 
is $L^p$ for any $p < +\infty$. 
By (5.15) $(\omega\wedge\omega)_{(4)} =
d*d\Phi(\omega\wedge\omega)_{(4)}$. 
 Hence 
\begin{eqnarray*}
&& \left(\omega\wedge\omega\right)_{(4)}
- d*d\left(u_4(\log \vert z \vert)z + \overline{u_4(\log \vert z
\vert)}\zbar\right)
\\
&=&d*d\left(\Phi(\omega\wedge\omega)_{(4)} 
- \left(u_4(\log \vert z \vert)z + \overline{u_4(\log \vert z
\vert)}\zbar\right)\right)
\end{eqnarray*}
is $L^p$ near $P_0$ for any $p < +\infty$.
From Corollary 5.3
$$
\omega_{(4)} -*d\left(u_4(\log \vert z \vert)z + \overline{u_4(\log \vert z
\vert)}\zbar\right)
$$
is H\"older continuous of class $C^{0+\alpha}$ for any $\alpha \in ]0, 1[$. 
Since $\deg u_4(s) \leq 1$, we have some constant $C_4 > 0$ such that 
$\left\vert\omega_{(4)}\left(\frac{d}{dz}\right)\right\vert \leq
C_4\vert\log\vert z\vert\vert$ near $P_0$, as was to be shown for the case
$m=4$.\par
Finally we suppose $m \geq 5$. From the inductive assumption we have 
some $u_{m-2}(s) \in H^{\otimes (m-2)}\otimes \bC[s]$ of degree 
$\leq (m-3)/2$ such that 
$$
\widetilde{\omega}_{({m-2})}:= 
\omega_{({m-2})} 
- *d\left(u_{m-2}(\log\vert z\vert)z+ \overline{u_{m-2}(\log\vert
z\vert)}\zbar\right)
$$
is H\"older continuous of class $C^{0+\alpha}$ 
for any $\alpha \in ]0, 1[$. 
If we choose a polynomial $\widehat{u}_{m-2}(s)$ satisfying 
$$
\widehat{u}'_{m-2}(s) + \frac12 \widehat{u}''_{m-2}(s) =
\frac{1}{4\pi\sqrt{-1}}u_{m-2}(s),
$$
then, from Lemma 5.5, we have 
$$
*d\left(u_{m-2}(\log\vert z\vert)z\right)\wedge
\left(\frac1{2\pi}d\theta\right)
= d*d \left(\widehat{u}_{m-2}(\log\vert z\vert)z \right)
$$
The value of the real form $\widetilde{\omega}_{({m-2})}$ at $P_0$ 
is given by 
$$
\widetilde{\omega}_{({m-2})}(P_0) 
= w_{m-2}(dz)_{P_0} +
\overline{w_{m-2}}(d\zbar)_{P_0}
$$
for some $w_{m-2} \in H^{\otimes m-2}\otimes \bC$. 
From the H\"older continuity of $\widetilde{\omega}_{({m-2})}$ 
$$
\vert z\vert^{-\alpha}\left(\widetilde{\omega}_{({m-2})} - w_{m-2}dz -
\overline{w_{m-2}}d\zbar\right)
$$
is bounded near $P_0$ for any $\alpha \in ]0, 1[$.
Hence 
$$
\left(\widetilde{\omega}_{({m-2})} - w_{m-2}dz -
\overline{w_{m-2}}d\zbar
\right)\wedge \left(\frac1{2\pi}d\theta\right)
$$
is $L^p$ near $P_0$ for any $p < +\infty$. 
If we define 
$$
u_m(s) := (\frac{s}{4\pi}w_{m-2}+ \widehat{u}_{m-2}(s)) I -
I (\frac{s}{4\pi}w_{m-2}+ \widehat{u}_{m-2}(s)),
$$
then 
$$
\omega_{(m-2)}\wedge\omega_{(2)} + \omega_{(2)}\wedge\omega_{(m-2)}
- d*d\left(u_m(\log \vert z \vert)z + \overline{u_m(\log \vert z
\vert)}\zbar\right)
$$
is $L^p$ near $P_0$ for any $p < +\infty$. 
From (5.16) and the inductive assumption for (4) 
the $2$-current $\sum_{p \neq 2, m-2} \omega_{(p)}\wedge\omega_{(m-p)}
= d*d\left(\sum_{p \neq 2, m-2} \omega_{(p)}\wedge\omega_{(m-p)}\right)$ 
is $L^p$ near $P_0$ for any $p < +\infty$. Therefore 
\begin{eqnarray*}
&& \left(\omega\wedge\omega\right)_{(m)}
- d*d\left(u_m(\log \vert z \vert)z + \overline{u_m(\log \vert z
\vert)}\zbar\right)\\
&=&d*d\left(\Phi(\omega\wedge\omega)_{(m)} - \left(u_m(\log \vert z
\vert)z +
\overline{u_m(\log \vert z \vert)}\zbar\right)\right)
\end{eqnarray*}
is also $L^p$ near $P_0$ for any $p < +\infty$. 
From Corollary 5.3 
$$
\omega_{(m)} -*d\left(u_m(\log \vert z \vert)z + \overline{u_m(\log \vert z
\vert)}\zbar\right)
$$
is H\"older continuous of class $C^{0+\alpha}$ for any $\alpha \in ]0, 1[$. 
Since $\deg u_m(s) = \deg \widehat{u}_{m-2}(s) \leq
\frac12(m-1)$,  we have some constant $C_m >0$ such that
$$
\left\vert\omega_{(m)}\left(\frac{d}{dz}\right)\right\vert \leq
C_m \left\vert \log\vert z\vert\right\vert^{[(m-1)/2]}.
$$ 
This completes the induction.
\end{proof}

\section{The Harmonic Magnus Expansion}

As in the previous sections let $(C, P_0, v)$ be a triple of 
a compact Riemann surface $C$ of genus $g \geq 1$, $P_0 \in C$ 
and $v \in T_{P_0}C\setminus\{0\}$. We define the map 
$$
\theta = \theta^{(C, P_0, v)}: \frak L(C, P_0, v) \to \T
$$
by the improper iterated integral
\begin{equation}
\theta(\gamma) := 1 + \int_\gamma \omega + \int_\gamma \omega\omega +
\int_\gamma \omega\omega\omega + \cdots =
1 + \sum^\infty_{m=1}\int_\gamma
\overbrace{\omega\omega\cdots \omega}^{\mbox{$m$ times}}
\label{6.1}\end{equation}
of the connection form $\omega$ constructed in \S5 along a loop 
$\ell \in \frak L(C, P_0, v)$.\par
We have to prove the improper integral converges. 
For any $\ell \in \frak L(C, P_0, v)$ and $1$-forms 
$\varphi_1, \dots, \varphi_m$ on $C \setminus \{P_0\}$ we have
\begin{eqnarray}
\left\vert\int_\ell\varphi_1 \cdots \varphi_m\right\vert
&\leq& \int_{1 > t_1\geq \cdots \geq t_m > 0}
\vert{\ell}^*\varphi_1\vert\times\cdots\times
 \vert{\ell}^*\varphi_m\vert \label{6.1}\\
&\leq& \left(\int^1_0\vert{\ell}^*\varphi_1\vert\right)\cdots
\left(\int^1_0\vert{\ell}^*\varphi_m\vert\right).\nonumber
\end{eqnarray}
Hence it suffices to show the following. 
As in the previous section we fix a complex coordinates $z$ 
of $C$ centered at $P_0$.
\begin{lemma}\label{Lemma 6.1} The $1$-form $\ell^*\omega_{(m)}$ is
integrable  on the interval $[0, 1]$
$$
\int^1_0 \vert\ell^*\omega_{(m)}\vert < +\infty
$$
for any $\ell \in \frak L(C, P_0, v)$ and $m \geq 1$.
\end{lemma}
\begin{proof} It is clear for the case $m = 1$. 
To study $\ell^*\omega_{(2)}$ 
we expand $x(t) = \Re z(\ell(t))$ and $y(t) = \Im z(\ell(t))$ 
in Taylor polynomials
\begin{eqnarray*}
&& x(t) = x_1t + x_2t^2 + o(t^2)\\
&& y(t) = y_1t + y_2t^2 + o(t^2).
\end{eqnarray*}
Then ${x_1}^2 + {y_1}^2 \neq 0$ because $({d\ell}/{dt})(0) = v \neq 0$. 
Substituting them into $d\theta = (xdy -ydx)/(x^2+y^2)$, we have 
$$
\ell^*d\theta = \frac{x_1y_2- x_2y_1+  o(1)}{{x_1}^2 + {y_1}^2 +o(1)}dt.
$$
This means $\ell^*d\theta$ is contiuous on $[0, 1[$. 
Similarly we find $\ell^*d\theta$ is contiuous on $]0, 1]$. 
Since $\widetilde{\omega}_{(2)} = 
\omega_{(2)} - \frac1{2\pi}Id\theta$ is of class $C^\infty$, 
$\ell^*\omega_{(2)}$ is continuous and integrable on the interval $[0, 1]$. 
\par
Now we study the case $m \geq 3$. 
From Proposition 5.4 (4) we have some constant $C'_m > 0$ such that 
\begin{eqnarray*}
\left\vert(\ell^*\omega_{(m)})\left(\frac{d}{dt}\right)\right\vert
&\leq& C'_m \vert\log t\vert^{[(m-1)/2]}\quad \text{near $t=0$, and}\\
&\leq& C'_m \vert\log (1-t)\vert^{[(m-1)/2]}\quad \text{near $t=1$.}
\end{eqnarray*}
Since $\int^1_0\vert\log t\vert^{[(m-1)/2]}dt < +\infty$,
the $1$-form 
$\ell^*\omega_{(m)}$ is integrable, as was to be shown. 
\end{proof}

Next we prove the homotopy invariance of $\theta(\ell)$. 
\begin{lemma}\label{Lemma 6.2}
If $\ell_0$ and $\ell_1 \in \frak L(C, P_0, v)$ are homotopic 
to each other in the sense of \S4, then we have $\theta(\ell_0) 
= \theta(\ell_1) \in \T$. 
\end{lemma}
\begin{proof} Let a piecewise $C^\infty$ map $L: [0, 1]\times [0, 1] \to 
C$ be a homotopy connecting $\ell_0$ and $\ell_1$ as in (4.2). 
For any $t \in [0, 1]$ we denote $k_t(s) := L(t, s)$, $s \in [0, 1]$. 
It suffices to show 
\begin{equation}
\lim_{t\downarrow 0}\left\vert\int_{k_t}\vert\omega_{(m)}\vert\right\vert
= \lim_{t\uparrow 1}\left\vert\int_{k_t}\vert\omega_{(m)}\vert\right\vert
= 0\label{6.3}
\end{equation}
for any $m \geq 1$. In fact, the connection form $\omega$ satisfies 
the integrability condition $d\omega = \omega \wedge\omega$ on 
$C\setminus \{P_0\}$ from (5.12). 
Hence we have 
\begin{gather*}
\left(1 +
\sum^\infty_{m=1}\int^{1-\epsilon}_\epsilon\overbrace
{({\ell_1}^*\omega)\cdots({\ell_1}^*\omega)}^m\right)
\left(1 +
\sum^\infty_{m=1}\int_{k_\epsilon}\overbrace
{\omega\cdots\omega}^m\right)\\
=
\left(1 +
\sum^\infty_{m=1}\int_{k_{1-\epsilon}}\overbrace
{\omega\cdots\omega}^{m}\right)\left(1+
\sum^\infty_{m=1}\int^{1-\epsilon}_\epsilon\overbrace
{({\ell_0}^*\omega)\cdots({\ell_0}^*\omega)}^{m}\right).
\end{gather*}
From (6.3), (6.2) and Lemma 6.1 we have $\theta(\ell_1)\cdot 1=
1\cdot\theta(\ell_0)$ as $\epsilon\downarrow 0$.\par
Now we will prove $\lim_{t\downarrow
0}\left\vert\int_{k_t}\vert\omega_{(m)}\vert\right\vert = 0$. One can
prove $\lim_{t\uparrow
1}\left\vert\int_{k_t}\vert\omega_{(m)}\vert\right\vert = 0$ in a similar
way. The $1$-forms $L^*\omegaone$ and $L^*\widetilde{\omega}_{(2)}$ are 
$C^\infty$ on $[0, 1]\times [0, 1]$ and the path $k_0$ is constant. 
Hence $\lim_{t\downarrow 0}\left\vert\int_{k_t}
\vert\omega_{(1)}\vert\right\vert = \lim_{t\downarrow 0}\left\vert\int_{k_t}
\vert\widetilde{\omega}_{(2)}\vert\right\vert = 0$. \par
From the second condition of (4.2) there exists a $C^\infty$ function
$\widehat{L}(t, s)$ satisfying $(z\circ L)(t, s) = t\widehat{L}(t, s)$ 
near $t=0$. 
By the third condition $\widehat{L}(0, s)$ is a nonzero constant. 
If we define $\widehat{k}_t: [0, 1] \to \bC$ by 
$\widehat{k}_t(s) := \widehat{L}(t, s)$, then 
${k_t}^*(d\theta) = {\widehat{k}_t}^*(d\theta)$ is of class $C^\infty$ on $[0,
1]\times [0, 1]$ and the path $\widehat{k}_0$ is constant. 
Hence $\lim_{t\downarrow 0}\left\vert\int_{k_t}
\vert d\theta\vert\right\vert = \lim_{t\downarrow 0}\left\vert\int_{k_t}
\vert\omega_{(2)}\vert\right\vert = 0$. 
\par
Consider the case $m \geq 3$. Clearly we have 
$(\partial (z\circ L)/\partial s)(t, s) = t (\partial \widehat{L}/\partial
s)(t, s)$. From Proposition 5.4 (4) we have some constant $C''_m > 0$ 
such that
$$
\left\vert({k_t}^*\omega_{(m)})\left(\frac{d}{ds}\right)\right\vert
\leq C''_m t\vert\log t\vert^{[(m-1)/2]}
$$
for any $(t, s) \in [0, 1]\times [0, 1]$. 
This implies $\lim_{t\downarrow 0}\left\vert\int_{k_t}\vert\omega_{(m)}\vert\right\vert
= 0$, as was to be shown.\qed
\end{proof}
Consequently the iterated integral (6.1) defines a group homomorphism
$$
\theta = \theta^{(C, P_0, v)}: \pi_1(C, P_0, v) \to 1+\T_1, \quad
[\ell] \mapsto \theta(\ell),
$$
which is a Magnus expansion of the free group $\pi_1(C, P_0, v)$ 
because of (5.8). We call it {\bf the harmonic Magnus expansion} 
of the triple $(C, P_0, v)$. \par

We conclude this section by studying how the harmonic Magnus expansion 
$\theta^{(C, P_0, v)}$ depends on the choice of the vector $v$. 
\begin{proposition}\label{Proposition 6.3} 
Let $\ell: [0, 1] \to C$ be a
piecewise 
$C^\infty$ path satisfying the conditions 
\begin{eqnarray*}
&& \ell(]0, 1[) \subset \{0 < \vert z\vert \ll 1 \},\\
&& \ell(0) = \ell(1) = P_0,\\
&& \frac{d\ell}{dt}(0) \neq 0, \quad\text{and}\quad \frac{d\ell}{dt}(1)
\neq 0.
\end{eqnarray*}
Then we have 
$$
\left(\theta(\ell) = \right) 1 +
\sum^\infty_{m=1}\int_\ell\overbrace{\omega\omega\cdots \omega}^m
 = \exp\left(\left(\frac1{2\pi}\int_\ell d\arg z\right)I\right) \in \T.
$$
In particular, the word $w_0$ given by a negative loop arround $P_0$ has 
its value
$$
\theta(w_0) = \exp(-I) \in \T.
$$
\end{proposition}
\begin{proof} For any $\rho \in ]0, \frac12[$ there exist some 
$a_\rho$ and $b_\rho \in \bC$ such that the path $t \in [\rho, 1-\rho]
\mapsto z\inv(e^{a_\rho t+b_\rho}) \in C$ is homotopic to the restriction 
$\ell\vert_{[\rho, 1-\rho]}$ relative to the endpoints $\{\rho, 1-\rho\}$. 
We define a piecewise $C^\infty$ path $\ell_\rho: [0, 1] \to C$ by 
$$
\ell_\rho := \left\{
\begin{array}{ll}
\ell, & \quad\text{on $[0, \rho] \cup [1-\rho, 1]$,}\\
z\inv(e^{a_\rho t+b_\rho}), & \quad\text{on $[\rho, 1-\rho]$.}
\end{array}
\right.
$$
From the homotopy invariance of $\theta^{(C, P_0, v)}$ we have 
$\theta(\ell) = \theta(\ell_\rho) = \lim_{\rho\downarrow0}\theta(\ell_\rho)$. 
Clearly we have $\frac1{2\pi}\int_{\ell_\rho}d\arg z =
\frac1{2\pi}\int_{\ell}d\arg z$, which we denote by $s_0$. \par
Since $\ell^*\omega_{(m)}$ is integrable on $[0, 1]$, we have 
$\lim_{\rho\downarrow0}\left(\int^\rho_0 +
\int^1_{1-\rho}\right)\vert\ell^*\omega_{(m)}\vert = 0$. 
Moreover we have some constants $C_1 > C_0 > 0$ such that 
${C_0} \rho\leq \vert z(\ell_\rho(t))\vert \leq {C_1}\rho$ for any 
$t \in [\rho, 1-\rho]$. If $m \geq 3$, then we have some 
constant $C'''_m > 0$ such that 
$$
\int^{1-\rho}_\rho\vert{\ell_\rho}^*\omega_{(m)}\vert \leq 
C'''_m\rho\vert\log\rho\vert^{[(m-1)/2]}
$$
from Proposition 5.4 (4).
Hence we have 
$$
\lim_{\rho\downarrow0}\int^1_0\vert{\ell_\rho}^*\omega_{(m)}\vert = 0
$$
for any $m \geq 3$. Clearly we have 
$$
\lim_{\rho\downarrow0}\left\vert\int^1_0\vert{\ell_\rho}^*\omega_{(1)}
\vert\right\vert
=\lim_{\rho\downarrow0}\left\vert\int^1_0\vert{\ell_\rho}^*\widetilde{
\omega}_{(2)}
\vert\right\vert = 0.
$$
Consequently, from (6.2), we obtain 
\begin{eqnarray*}
&& \lim_{\rho\downarrow0}\theta(\ell_\rho) 
 = \lim_{\rho\downarrow0}\sum_{m \geq 0 \atop \text{$m$: even}}
\int_{\ell_\rho}\overbrace{\left(\frac1{2\pi}Id\arg z\right)
\cdots\left(\frac1{2\pi}Id\arg z\right)}^{\text{$m/2$
times}}\\
&=& \sum_{m \geq 0 \atop \text{$m$: even}}I^{\otimes (m/2)}
\int^{s_0}_{0}\overbrace{(dt)\cdots(dt)}^{\text{$m/2$
times}}
= \sum_{m \geq 0 \atop \text{$m$: even}}
\frac1{(m/2)!}{s_0}^{m/2}I^{\otimes (m/2)}= \exp(s_0I),
\end{eqnarray*}
as was to be shown.\end{proof}

\section{Quasiconformal Variation of Harmonic Magnus Expansions}

Let $(C, P_0, v)$ be a triple of genus $g \geq 1$, 
and $z$ a complex coordinate of $C$ centered at $P_0$
as in the previous sections. 
Moreover let $(C_t, {P_0}^t, v_t)$, $t \in \bR$, $\vert t \vert \ll 1$, 
be a $C^\infty$ family of triples of genus $g$ satisfying 
$(C_t, {P_0}^t, v_t)\bigr\vert_{t=0} = (C, {P_0}, v)$. 
Then the family $(C_t, {P_0}^t, v_t)$ is trivial as a $C^\infty$ 
fiber bundle over an interval near $t = 0$. 
Hence we have a $C^\infty$ family of $C^\infty$ diffeomorphisms 
$$
f^t: (C, P_0, v) \to (C_t, {P_0}^t, v_t)
$$
satisfying $f^0 = 1_{(C, P_0, v) }$. 
We may assume $(df^t)_{P_0}: (T_{\bR}C)_{P_0} \to (T_{\bR}C_t)_{P^t_0}$ 
preserves the almost complex structures on the tangent spaces at 
the points $P_0$ and $P^t_0$ for any $t$. This means 
\begin{equation}
{f^t}_\zbar(P_0) = 0 
\label{7.1}
\end{equation}
for any $t \in \bR$, $\vert t \vert \ll 1$. \par
In this section we compute the first quasiconformal variation
$$
\thetadot := \frac{d}{dt}\Bigr\vert_{t=0}\theta_t:
\pi_1(C, P_0, v) \to \T
$$
of the $C^\infty$ family of Magnus expansions
$$
\theta_t := (f^t)^*\left(\theta^{(C_t, P^t_0, v_t)}\right): 
\pi_1(C, P_0, v) \to \T.
$$
In general, if $\bigcirc = \bigcirc_t$ is a ``function" in 
$t \in \bR$, $\vert t \vert \ll 1$, then we write simply 
$$
\overset\centerdot\bigcirc := \frac{d}{dt}\Bigr\vert_{t=0}\bigcirc_t.
$$
For example, we denote 
$$
\mudot := \frac{d}{dt}\Bigr\vert_{t=0}\mu(f^t).
$$
Here $\mu(f^t)$ is the complex dilatation of the diffeomorphism $f^t$. 
Let $z_1$ be a complex coordinate of $C$, and $\zeta_1$ of $C_t$. 
The complex dilatation $\mu(f^t)$ is defined locally by 
$$
\mu(f^t) = \mu(f^t)(z_1)\frac{d}{dz_1}\otimes d\overline{z_1} = 
\frac{(\zeta_1\circ
f^t)_{\overline{z_1}}}{(\zeta_1\circ
f^t)_{z_1}}\frac{d}{dz_1}\otimes d\overline{z_1},
$$
which does not depend on the choice of the coordinates $z_1$ and $\zeta_1$. 
From (7.1) we have 
\begin{equation}
\mudot(P_0) = 0. 
\label{7.2}
\end{equation}
To state our result we introduce linear maps $N$ and 
$\Check N: \T_1 \to \T_1$ by 
$$
N\vert_{H^{\otimes m}} :=
1 + \varepsilon + \varepsilon^2 + \cdots +
\varepsilon^{m-1} = {\sum}^{m-1}_{k=0}\,\varepsilon^k, \quad\text{and}\quad
\Check N\vert_{H^{\otimes m}} := \frac1m N,
$$
respectively.\par
Applying the map $N$ to the $\T_2$-part of the square 
$\omega'\omega' \in C^\infty(C \setminus \{P_0\};
(T^*C)^{\otimes 2})\wotimes \T_2$ of the $(1, 0)$-part of the
connection form $\omega$, we obtain a covariant tensor 
$N(\omega'\omega') \in C^\infty(C \setminus \{P_0\};
(T^*C)^{\otimes 2})\wotimes \T_2$.
 
\begin{lemma}\label{Lemma 7.1} The covariant tensor $N(\omega'\omega')$
is  holomorphic on $C \setminus \{P_0\}$. 
Except for the case $m=4$ $N(\omega'\omega')_{(m)}$ has a
pole of order $\leq 1$ at the point $P_0$, while 
$N(\omega'\omega')_{(4)}$ a
pole of order $\leq 2$ at $P_0$.
\end{lemma}
\begin{proof} On the complement $C \setminus \{P_0\}$ we have $d\omega =
\omega\wedge\omega$, i.e., 
$$
\overline{\partial}\omega'_{(m)}
= \frac12 \left(\omega'\wedge\omega'' 
+ \omega''\wedge\omega'\right)_{(m)}.
$$
This implies
\begin{eqnarray*}
&& 2 \overline{\partial} N(\omega'\omega')_{(m)}
= 2\overline{\partial}
N\sum_{p+q = m}\omega'_{(p)}\omega'_{(q)}\\
&=& N\sum_{p_1+p_2+q = m}
\left(\omega'_{(p_1)}\wedge\omega''_{(p_2)}\right)\omega'_{(q)} +
\left(\omega''_{(p_1)}\wedge\omega'_{(p_2)}\right)\omega'_{(q)}
\\
&& + N\sum_{p+q_1+q_2 = m}\omega'_{(p)}
\left(\omega'_{(q_1)}\wedge\omega''_{(q_2)}\right) +
\omega'_{(p)}\left(\omega''_{(q_1)}\wedge\omega'_{(q_2)}\right)
\\
&=& N\sum_{p_1+p_2+q = m}
\varepsilon^{-q}\left(\omega'_{(q)}\omega'_{(p_1)}\otimes\omega''_{(p_2)}\right)-
\varepsilon^{p_1}\left(\omega'_{(p_2)}\omega'_{(q)}\otimes\omega''_{(p_1)}\right)
\\
&& + N\sum_{p+q_1+q_2 = m}
\omega'_{(p)}\omega'_{(q_1)}
\otimes\omega''_{(q_2)}
-\varepsilon^{-q_2}\left(\omega'_{(q_2)}\omega'_{(p)}\otimes\omega''_{(q_1)}\right)
\\
&=& N\sum_{a+b+c=m}\left(\varepsilon^{-a}- \varepsilon^c +1
-\varepsilon^{-a}\right)
\omega'_{(a)}\omega'_{(b)}\otimes\omega''_{(c)}
= 0
\end{eqnarray*} 
Hence $N(\omega'\omega')$ is holomorphic on $C \setminus \{P_0\}$.\par
If $m \neq 4$, then, from Proposition 5.4, we have 
$\lim_{z\to 0}\vert z\vert^{2}N(\omega'\omega')_{(m)} = 0$, and 
$\lim_{z\to 0}\vert z\vert^{2+\epsilon}N(\omega'\omega')_{(4)} = 0$
for any $\epsilon > 0$, which implies the second half of the asserion. 
This proves the lemma.
\end{proof}

Our result is 

\begin{theorem}\label{Theorem 7.2}
$$
\thetadot(\gamma) = \inter\left(\int_C
2\Re\left(\left(N(\omega'\omega')- 2\omega'_{(1)}\omega'_{(1)}\right)\mudot\right)\right)
\theta(\gamma) \in \T
$$
for any $\gamma \in \pi_1(C, P_0, v)$.
\end{theorem}

Here we regard $N(\omega'\omega')\mudot$ as an integrable $2$-current on 
$C$ by (7.2) and Lemma 7.1. The second homogeneous 
term of $N(\omega'\omega')$
is 
 $2{\omega'}_{(1)}{\omega'}_{(1)}$, which coincides with the first
quasiconformal variation of the period matrix of Riemann surfaces 
given by Rauch's variational formula \cite{R}. \par
The rest of this section is devoted to the proof of the theorem. 
We begin our computation by recalling the first variation of 
the Hodge $*$-operators. We write locally 
$$
\mudot = \mudot(z_1)\frac{d}{dz_1}\otimes d\overline{z_1}
$$
in a complex coodinate $z_1$ of $C$.
\begin{lemma}\label{Lemma 7.3} 
$$
\frac{d}{dt}\Bigr\vert_{t=0}(f^t)^**_{C_t}((f^t)^{-1})^*dz_1 =
-2\sqrt{-1}\mudot(z_1) d\overline{z_1}
$$
\end{lemma}
\begin{proof} Let $\zeta^t_1$ be a complex coordinates of $C_t$ which 
is of class $C^\infty$ in the variable $t$, and satisfies 
$\zeta^0_1 = z_1$. Then we have 
\begin{eqnarray*}
&& (f^t)^**_{C_t}((f^t)^{-1})^*dz_1 - *dz_1\\
&=& (f^t)^*(*_{C_t} + \sqrt{-1})((f^t)^{-1})^*dz_1\\
&=& (f^t)^*(*_{C_t} + \sqrt{-1})\left(
(z_1\circ(f^t)^{-1})_{\zeta^t_1} d\zeta^t_1 + 
(z_1\circ(f^t)^{-1})_{\overline{\zeta^t_1}} d\overline{\zeta^t_1}
\right)\\
&=& 2\sqrt{-1}(f^t)^*\left(
(z_1\circ(f^t)^{-1})_{\overline{\zeta^t_1}} d\overline{\zeta^t_1}
\right)\\
&=&
2\sqrt{-1}\left((z_1\circ(f^t)^{-1})_{\overline{\zeta^t_1}}\circ
f^t\right)
\left((\overline{\zeta^t_1\circ f^t})_{z_1}dz_1 + (\overline{\zeta^t_1\circ
f^t})_{\overline{z_1}}d\overline{z_1}
\right).
\end{eqnarray*}
Since $\lim_{t\to 0}(z_1\circ(f^t)^{-1})_{\overline{\zeta^t_1}} = 
\lim_{t\to 0}(\overline{\zeta^t_1\circ f^t})_{z_1} = 0$, we obtain 
$$
\lim_{t\to 0}\frac1t\left((f^t)^**_{C_t}((f^t)^{-1})^*dz_1 - *dz_1\right) 
= 2\sqrt{-1}\left(\lim_{t\to 0}\frac1t\left((z_1\circ(f^t)^{-1})_{\overline{\zeta^t_1}}\circ
f^t\right)\right)d\overline{z_1}.
$$
Hence it suffices to show
$$
\lim_{t\to 0}\frac1t\left((z_1\circ(f^t)^{-1})_{\overline{\zeta^t_1}}\circ
f^t\right) = -\mudot(z_1).
$$
Now we have 
\begin{eqnarray*}
0 &=& (z_1)_{\overline{z_1}} =
(z_1\circ(f^t)^{-1}\circ f^t)_{\overline{z_1}}\\
&=& \left((z_1\circ(f^t)^{-1})_{\zeta^t_1}\circ f^t\right)(\zeta^t_1\circ
f^t)_{\overline{z_1}}  +
\left((z_1\circ(f^t)^{-1})_{\overline{\zeta^t_1}}\circ f^t\right)
(\overline{\zeta^t_1\circ f^t})_{\overline{z_1}}.
\end{eqnarray*}
This implies
\begin{eqnarray*}
&& \lim_{t\to
0}\frac1t\left((z_1\circ(f^t)^{-1})_{\overline{\zeta^t_1}}\circ
f^t\right) =\lim_{t\to
0}\frac1t\left((z_1\circ(f^t)^{-1})_{\overline{\zeta^t_1}}\circ
f^t\right)(\overline{\zeta^t_1\circ f^t})_{\overline{z_1}} \\ 
&=& -
\lim_{t\to 0}\frac1t\left((z_1\circ(f^t)^{-1})_{\zeta^t_1}\circ
f^t\right) (\zeta^t_1\circ f^t)_{\overline{z_1}} =  - \lim_{t\to
0}\frac1t (\zeta^t_1\circ f^t)_{\overline{z_1}}\\
&=& - \lim_{t\to
0}\frac1t \mu(f^t)(z_1)(\zeta^t_1\circ f^t)_{z_1} = -\mudot(z_1),
\end{eqnarray*}
as was to be shown.
\end{proof}

We define a linear operator $S = S[\mudot]: A^1(C) \to A^1(C)$ by 
$$
S(\varphi) = S(\varphi') + S(\varphi'') :=
-2\varphi'\mudot - 2\varphi''\overline{\mudot},
$$
for $\varphi = \varphi' + \varphi''$,  $\varphi' \in A^{1, 0}(C)$, 
$\varphi'' \in A^{0, 1}(C)$. 
Here, if $\varphi' = \varphi'(z_1)dz_1$ and $\varphi'' =
\varphi''(z_1)d\overline{z_1}$ locally in a complex coordinate $z_1$ of
$C$, then we have locally
$$
-2\varphi'\mudot = -2\varphi'(z_1)\mudot(z_1)d\overline{z_1},
\quad\text{and}\quad -2\varphi''\overline{\mudot} = 
-2\varphi''(z_1)\overline{\mudot(z_1)}dz_1.
$$
Clearly $S(A^{1, 0}(C)) \subset A^{0, 1}(C)$,  
$S(A^{0, 1}(C)) \subset A^{1, 0}(C)$, and 
$S$ is a real operator, i.e.,
$\overline{S} = S$. From Lemma 7.3 we have 
\begin{equation}
\stardot = *S = -S*: A^1(C) \to A^1(C).
\label{7.3}
\end{equation}
It is easy to show 
\begin{equation}
(S\varphi)\wedge\psi = - \varphi\wedge(S\psi) 
\label{7.4} 
\end{equation}
for any $\varphi$ and $\psi \in A^1(C)$. 
From (7.2) $\omega'\mudot$ and $\omega''\overline{\mudot}$ are bounded 
near the point $P_0$. Hence we have 
$$
S\omega \in E^1(C)\wotimes\T,
$$
where $E^1(C)$ is the function space introduced in (5.4). \par
To describe the first variation of the connection form $\omega$ 
we introduce a linear map $U: E^1(C)\wotimes\T \to A^1(C)
\wotimes\T$ by 
$$
U\varphi := [\omega, \Phi d*\varphi] = \omega (\Phi d*\varphi) - (\Phi
d*\varphi)\omega
$$
for $\varphi \in E^1(C)\wotimes\T$.
This $U$ appears only in this section, and is completely different 
from the $Sp_{2g}(\bR)$-module $\La^3H/H$ in \S8.

\begin{lemma}\label{Lemma 7.4} 
For any $\varphi \in E^1(C)\wotimes\T$ we have
$$
U\varphi \in E^1(C)\wotimes\T = \left({\bigcap}_{2 < p <
\infty}L^p\left(C;
\, (T^*_\bR C)\otimes\bC\right)\right)\wotimes\T,\leqno{(1)}
$$
$$
d(U\varphi) = (U\varphi - d\Phi d*\varphi)\wedge\omega +
\omega\wedge(U\varphi - d\Phi d*\varphi).\leqno{(2)}
$$
\end{lemma}

\begin{proof} (1) From Theorem 5.1 we have $\Phi d*\varphi \in E^0_0(C,
P_0)$.  Clearly $U\varphi \in E^1(C- \{P_0\})\wotimes\T$. It suffices to
show that 
$\Phi d*\varphi$ is $L^p$ near $P_0$ for any $p\in ]2, \infty[$. \par
Choose $p' \in ]p, \infty[$. Then $\vert z\vert^{-1 + 2/p'}\vert\Phi
d*\varphi\vert$ and $\vert z\vert\vert\omega\vert$ are bounded near $P_0$. 
Now we have 
$$
\int_{\vert z\vert \leq 1}\left(\vert z\vert^{-2/p'}\right)^{p}\vert
dz\wedge d\zbar\vert = 4\pi\int^1_0 r^{1-2{p}/p'}dr = 
4\pi\left[r^{2-2{p}/p'}\right]^{r=1}_{r=0} < +\infty.
$$
Hence $\omega (\Phi d*\varphi)$ and $(\Phi d*\varphi)\omega$ are 
$L^p$ near $P_0$, and so is $U\varphi$, as was to be shown.\par
(2) From $(\Phi d*\varphi)(P_0) = 0$ we have $(\Phi d*\varphi)\delta_0 = 0
\in A^2(C)$. Hence 
\begin{eqnarray*}
 d(U\varphi) 
&=& \omega\wedge\omega(\Phi d*\varphi) -
\omega\wedge(d\Phi d*\varphi) -
(d\Phi d*\varphi)\wedge\omega
- (\Phi d*\varphi)\omega\wedge\omega
\\
&=& (U\varphi)\wedge\omega + \omega\wedge(U\varphi) - \omega\wedge(d\Phi
d*\varphi) -
(d\Phi d*\varphi)\wedge\omega
\\
&=& (U\varphi - d\Phi d*\varphi)\wedge\omega +\omega\wedge(U\varphi -
d\Phi d*\varphi).
\end{eqnarray*}
\end{proof}

Since $U\left(E^1(C)\wotimes\T_p\right) \subset 
E^1(C)\wotimes\T_{p+1}$
for any $p \geq 0$, we may consider the operator 
$$
(1-U)\inv: E^1(C)\wotimes\T \to E^1(C)\wotimes\T
$$
and so the tensor $(1-U)\inv(-S\omega) \in E^1(C)\wotimes\T$. 

\begin{lemma}\label{Lemma 7.5}
$$
(1-\varepsilon)\int_C(1-U)\inv(-S\omega)\wedge\omegaone = 0. 
$$
\end{lemma}
\begin{proof} 
From (7.4) we have $\int_C(S\omega)\wedge\omega + \int_C\omega\wedge(S\omega)
= 0$. Using Lemma 7.4 (2) for $\varphi = (1-U)\inv(-S\omega)$ one computes
\begin{eqnarray*}
0 &=& \int_C d(U(1-U)\inv(-S\omega))\\
&=& \int_C (U(1-U)\inv(-S\omega) - 
d\Phi d*(1-U)\inv(-S\omega))\wedge\omega\\
&& + \int_C \omega\wedge(U(1-U)\inv(-S\omega) - 
d\Phi d*(1-U)\inv(-S\omega))\\
&=& \int_C (S\omega + (1-U)\inv(-S\omega) - 
d\Phi d*(1-U)\inv(-S\omega))\wedge\omega\\
&& + \int_C \omega\wedge(S\omega + (1-U)\inv(-S\omega) - 
d\Phi d*(1-U)\inv(-S\omega))\\
&=& \int_C (\harmonic + *d\Phi d)(1-U)\inv(-S\omega)\wedge\omega
 + \int_C \omega\wedge(\harmonic + *d\Phi d)(1-U)\inv(-S\omega)\\
&=& \int_C (1-U)\inv(-S\omega)\wedge(\harmonic + d\Phi d*)\omega
 + \int_C (\harmonic + d\Phi d*)\omega\wedge(1-U)\inv(-S\omega)\\
&=& \int_C (1-U)\inv(-S\omega)\wedge\omegaone
 + \int_C \omegaone\wedge(1-U)\inv(-S\omega)\\
&=& (1-\varepsilon)\int_C (1-U)\inv(-S\omega)\wedge\omegaone.
\end{eqnarray*}
\end{proof}

We define $W \in E^1(C)\wotimes \T$ by 
$$
W := (1-U)\inv(-S\omega) +S\omegaone
= (1-U)\inv(-S\omega) - ((1-U)\inv(-S\omega) )_{(1)}.
$$
From Lemma 7.5 and (4.5) we have 
\begin{equation}
\inter\left(\int_CW\wedge\omegaone\right)I = 0.
\label{7.5}
\end{equation}

Now we can compute the first variation of the connection form $\omega$. 
\begin{theorem}\label{Theorem 7.6}
$$
\omegadot = -\inter\left(\int_CW\wedge\omegaone\right)\omega
- U(1-U)\inv(-S\omega) + d\Phi d*(1-U)\inv(-S\omega).
$$
\end{theorem}

First we compute the $d$-exact and the harmonic parts of the 
variation $\omegadot$. 
\begin{lemma}\label{Lemma 7.7}
$$
d\Phi d*\omegadot = d\Phi d*(-S\omega),\leqno{(1)}
$$
$$
\harmonic\omegadot = -\harmonic S(\omega -\omegaone).\leqno{(2)}
$$
\end{lemma}
\begin{proof} (1) Differentiating $d*\omega = 0$, we have
$0 = d\stardot\omega + d*\omegadot = d*S\omega + d*\omegadot$. 
Hence $d\Phi d*\omegadot = d\Phi d*(-S\omega)$, as was to be shown.\par
(2) Consider the pullback $\omega^t_{(1)} :=
(f^t)^*\left(\omegaone^{C_t}\right)$ of the harmonic form $\omegaone^{C_t}$ on
$C_t$.  The difference $\omega^t_{(1)} - \omegaone$ is $d$-exact, 
which implies $\omega^t_{(1)} - \omegaone = 
d\Phi d*(\omega^t_{(1)} - \omegaone) = d\Phi d*(\omega^t_{(1)})$. 
Hence, from (1), we obtain
\begin{equation}
\omegadot_{(1)} = d\Phi d*\omegadot_{(1)} =
d\Phi d*(-S\omegaone). 
\label{7.6}
\end{equation}
Differentiating $I = \int_C\omega\wedge\omegaone$, we obtain 
$0 = \int_C\omegadot\wedge\omegaone + \int_C\omega\wedge\omegadot_{(1)}$. 
Hence, by (7.6), we obtain 
\begin{eqnarray*}
\harmonic\omegadot 
&=& -\left(\int_C\omegadot\wedge\omegaone\right)\cdot\omegaone
= \left(\int_C\omega\wedge\omegadot_{(1)}\right)\cdot\omegaone\\
&=& \left(\int_C\omega\wedge d\Phi d*(-S\omegaone)\right)\cdot\omegaone
 = \left(\int_C(\omega-\omegaone)\wedge (-S\omegaone)\right)\cdot\omegaone\\
&=& \left(\int_CS(\omega-\omegaone)\wedge \omegaone\right)\cdot\omegaone
= -\harmonic S(\omega-\omegaone).\\
\end{eqnarray*}
This completes the proof of the lemma.
\end{proof}

\begin{proof}[Proof of Theorem 7.6]
We study a $1$-form $\Theta \in
A^1(C)\wotimes\T$ defined by 
$$
\Theta := \omegadot + \inter\left(\int_CW\wedge \omegaone\right)\omega 
+ U(1-U)\inv(-S\omega)
$$
instead of $\omegadot$ itself. 
Then we have 
\begin{eqnarray}
& \harmonic \Theta = 0, \label{7.7}\\
& d\Phi d* \Theta = d\Phi d*(1-U)\inv(-S\omega). \label{7.8}
\end{eqnarray}
In fact, from Lemma 7.7(2), we have 
\begin{eqnarray*}
\harmonic\Theta &=& -\harmonic S(\omega-\omegaone) 
+ \inter\left(\int_CW\wedge \omegaone\right)\omegaone 
+ \harmonic\left( U(1-U)\inv(-S\omega)\right)\\ 
&=& \harmonic( S\omegaone) 
- \harmonic W
+ \harmonic\left( (1-U)\inv(-S\omega)\right) = \harmonic (W-W) = 0.
\end{eqnarray*}
Moreover, from Lemma 7.7 (1), we have 
$d\Phi d*\Theta = d\Phi d*(-S\omega) + d\Phi d* U(1-U)\inv(-S\omega)= d\Phi d*
(1-U)\inv(-S\omega)$, 
as was to be shown. \par
Now, from (7.5), we have $\inter\left(\int_CW\wedge\omegaone\right)d\omega 
= \left(\inter\left(\int_CW\wedge\omegaone\right)\omega\right)\wedge\omega +
\omega\wedge\left(\inter\left(\int_CW\wedge\omegaone\right)\omega\right)$. 
Since $f^t(P_0) = P^t_0$, we have $d\omegadot = \omegadot\wedge\omega +
\omega\wedge\omegadot$. Hence, by Lemma 7.4 (2), we have 
\begin{eqnarray*}
&& d\Theta =
d\left(\omegadot +
\inter\left(\int_CW\wedge \omegaone\right)\omega 
+ U(1-U)\inv(-S\omega)\right)
\\
&=& \omegadot\wedge\omega + \omega\wedge\omegadot 
+ \left(\inter\left(\int_CW\wedge \omegaone\right)\omega\right)\wedge\omega 
+ \omega\wedge\left(\inter\left(\int_CW\wedge \omegaone\right)\omega\right)\\
&& + \left(U(1-U)\inv(-S\omega) -
d\Phi d*(1-U)\inv(-S\omega)\right)\wedge\omega\\
&& + \omega\wedge\left(U(1-U)\inv(-S\omega) -
d\Phi d*(1-U)\inv(-S\omega)\right)
\\
&=& \left(\Theta - d\Phi d*(1-U)\inv(-S\omega)\right)\wedge\omega
 + \omega\wedge\left(\Theta - d\Phi d*(1-U)\inv(-S\omega)\right).
\end{eqnarray*}
From (7.7) and (7.8) follows
$\Theta - d\Phi d*(1-U)\inv(-S\omega) = *d\Phi d\Theta$.
Hence we obtain
$$
d\Theta = (*d\Phi d\Theta)\wedge\omega + \omega\wedge(*d\Phi d\Theta).
$$
This proves $(d\Theta)_{(m)} = 0$ inductively on $m \geq 1$, 
since $(d\Theta)_{(0)} = 0$. Consequently we have 
\begin{eqnarray*}
&& \omegadot + \inter\left(\int_CW\wedge
\omegaone\right)\omega  + U(1-U)\inv(-S\omega)\\
&=& \Theta = d\Phi d*\Theta = d\Phi d*(1-U)\inv(-S\omega),
\end{eqnarray*}
which completes the proof of the theorem.\end{proof}

The theorem, together with the formulae (1.5), implies

\begin{proposition}\label{Proposition 7.8}
$$
\thetadot(\gamma) = - \inter\left(\int_CW\wedge
\omegaone\right)\theta(\gamma)
$$
for any $\gamma \in \pi_1(C, P_0, v)$.
\end{proposition}

\begin{proof}
Since $\Phi d*(1-U)\inv(-S\omega)(P_0) = 0$, the variation
$$
\left(\int_\gamma
\overbrace{\omega\omega\cdots \omega}^{m}\right)^\centerdot =
\sum^m_{i=1}\int_\gamma\omega\cdots\overset{i}{\breve{\omegadot}}\cdots
\omega
$$
is equal to 
\begin{eqnarray*}
&& - \sum^m_{i=1}\int_\gamma
\omega\cdots\left(\inter\left(\int_CW\wedge\omegaone\right)\omega\right)
\cdots\omega\\
&& - \sum^m_{i=1}\int_\gamma
\omega\cdots U(1-U)\inv(-S\omega)\cdots\omega\\
&& + \sum^m_{i=1}\int_\gamma
\omega\cdots d\Phi d*(1-U)\inv(-S\omega)\cdots\omega\\
&=& -\inter\left(\int_CW\wedge\omegaone\right) \left(\int_\gamma
\overbrace{\omega\omega\cdots \omega}^{m}\right)
 - \sum^m_{i=1}\int_\gamma
\omega\cdots U(1-U)\inv(-S\omega)\cdots\omega\\
&& + \sum^m_{i=2}\int_\gamma
\overbrace{\omega\cdots\omega}^{i-2} \left(\omega\Phi
d*(1-U)\inv(-S\omega)\right)\overbrace{\omega\cdots\omega}^{m-i}\\
&& - \sum^{m-1}_{i=1}\int_\gamma
\overbrace{\omega\cdots\omega}^{i-1} \left(\Phi
d*(1-U)\inv(-S\omega)\omega\right)\overbrace{\omega\cdots\omega}^{m-i-1}\\
&=& -\inter\left(\int_CW\wedge\omegaone\right) \left(\int_\gamma
\overbrace{\omega\omega\cdots \omega}^{m}\right)
 - \sum^m_{i=1}\int_\gamma
\omega\cdots U(1-U)\inv(-S\omega)\cdots\omega\\
&& + \sum^{m-1}_{i=1}\int_\gamma
\overbrace{\omega\cdots\omega}^{i-1}
\left(U(1-U)\inv(-S\omega)\right)\overbrace
{\omega\cdots\omega}^{m-i-1},
\end{eqnarray*}
from the formulae (1.5). 
Consequently we obtain 
$$
\thetadot(\gamma) = \sum^\infty_{m=1}\left(\int_\gamma
\overbrace{\omega\omega\cdots \omega}^{m}\right)^\centerdot
= - \inter\left(\int_CR\wedge
\omegaone\right)\theta(\gamma).
$$
This proves the proposition.
\end{proof}

In order to identify the integral $\int_CW\wedge
\omegaone$ we introduce some additional notations. 
We define the $1$-current $\widehat{\omega} \in A^1(C)\wotimes\T$ 
by 
$$
\widehat{\omega}_{(m)} := m\omega_{(m)},
$$
and the operator $V: A^1(C)\wotimes \T \to A^1(C)\wotimes \T$ by
$$
V(\varphi) := *d\Phi(\omega\wedge\varphi+\varphi\wedge\omega)
$$ 
for $\varphi \in A^1(C)\wotimes\T$. Then we have

\begin{lemma}\label{Lemma 7.9}
$$
(\widehat{\omega}\wedge\omega + \omega\wedge \widehat{\omega})_{(m)} 
= m(\omega\wedge\omega)_{(m)}\leqno{(1)}
$$
$$
\check N(\widehat{\omega}'\omega' + \omega'\widehat{\omega}') =
N(\omega'\omega')\leqno{(2)}
$$
$$
\omegaone = (1-V) \widehat{\omega}\leqno{(3)}
$$
\end{lemma}

\begin{proof} (1) From the definition of $\widehat{\omega}$ we have 
\begin{eqnarray*}
&& (\widehat{\omega}\wedge\omega + \omega\wedge \widehat{\omega})_{(m)} =
{\sum}^{m-1}_{p=1}\widehat{\omega}_{(p)}\wedge\omega_{(m-p)}
+ \omega_{(p)}\wedge\widehat{\omega}_{(m-p)}\\
&=& {\sum}^{m-1}_{p=1}(p + m-p)\omega_{(p)}\wedge\omega_{(m-p)}
= m(\omega\wedge\omega)_{(m)}.
\end{eqnarray*}
\par
(2) can be proved in a similar way to (1).\par
(3) From (1) we have 
\begin{eqnarray*}
&& \left((1-V)\widehat{\omega}\right)_{(m)}
= \left(\widehat{\omega} - *d\Phi(\widehat{\omega}\wedge\omega +
\omega\wedge\widehat{\omega})\right)_{(m)} \\
&=& m\omega_{(m)} - m*d\Phi(\omega\wedge\omega)_{(m)} =
 m\omega_{(m)} - m(\omega - \omegaone)_{(m)}
 = \delta_{1m}\omegaone.
\end{eqnarray*}
\end{proof}

Now we can identify the integral $\int_C(1-U)\inv(-S\omega)\wedge\omegaone = 
\int_CW\wedge\omegaone - \int_C(S\omegaone)\wedge\omegaone$. 
\begin{lemma}\label{Lemma 7.10}
$$
\int_C(1-U)\inv(-S\omega)\wedge\omegaone =
-2\int_C\Re\left(N(\omega'\omega')\mudot\right) \in \T.
$$
\end{lemma}
\begin{proof}
First we remark $V\widehat{\omega}_{(m)} = m(\omega- \omegaone)_{(m)}$ is 
$L^q$ for any $q \in ]0, 2[$ and $(1-U)\inv(-S\omega)$ is $L^p$ for any 
$p \in ]2, \infty[$. Hence the product $(1-U)\inv(-S\omega)\wedge V
 \widehat{\omega}$ is integrable on $C$. Since $\Phi
d*(1-U)\inv(-S\omega)(P_0) = 0$, we have 
\begin{eqnarray*}
&& \check N \int_C(1-U)\inv(-S\omega)\wedge V
 \widehat{\omega}\\
&=& \check N
\int_C(1-U)\inv(-S\omega)\wedge*d\Phi(\widehat{\omega}\wedge\omega
+\omega\wedge\widehat{\omega})
\\
&=& \check N \int_Cd\Phi
d*(1-U)\inv(-S\omega)\wedge*d\Phi(\widehat{\omega}\wedge\omega
+\omega\wedge\widehat{\omega})
\\
&=& -\check N \int_C\Phi d*(1-U)\inv(-S\omega)\wedge
d*d\Phi(\widehat{\omega}\wedge\omega +\omega\wedge\widehat{\omega})
\\
&=& -\check N \int_C\Phi d*(1-U)\inv(-S\omega)\wedge
(\widehat{\omega}\wedge\omega +\omega\wedge\widehat{\omega})
\\
&=& \check N \int_C\left(\omega\Phi d*(1-U)\inv(-S\omega)- \Phi
d*(1-U)\inv(-S\omega)\omega\right)\wedge \widehat{\omega}
\\
&=& \check N \int_CU(1-U)\inv(-S\omega)\wedge
\widehat{\omega}.
\end{eqnarray*}
From Lemma 7.5 we have
$$
\int_C(1-U)\inv(-S\omega)\wedge\omegaone = \check
N \int_C(1-U)\inv(-S\omega)\wedge\omegaone.
$$
From Lemma 7.9 (2), (3) and what we have shown above, the RHS is equal to 
\begin{eqnarray*}
&& \check N\int_C(1-U)\inv(-S\omega)\wedge(1-V)\,
 \widehat{\omega}\\
&=& \check N\int_C(1-U)\inv(-S\omega)\wedge
 \widehat{\omega}
 - \check N\int_C(1-U)\inv(-S\omega)\wedge
V \widehat{\omega}
 \\
&=& \check N\int_C(1-U)\inv(-S\omega)\wedge
 \widehat{\omega}
 - \check N\int_CU(1-U)\inv(-S\omega)\wedge
 \widehat{\omega}
 \\
&=& \check N \int_C
(-S\omega)
\wedge\widehat{\omega}
=  \check N \int_C
2(\omega'\mudot + \omega''\overline{\mudot})
\wedge(\widehat{\omega}' + \widehat{\omega}'')\\
&=& -2\check N \int_C (\omega'\widehat{\omega}'\mudot
+ \omega''\widehat{\omega}''\overline{\mudot})
= -4\check N \int_C \Re(\omega'\widehat{\omega}'\mudot)\\
&=& -2 \int_C \Re(\check
N(\omega'\widehat{\omega}'+\widehat{\omega}'\omega')\mudot)
= -2\int_C\Re(N(\omega'\omega')\mudot),
\end{eqnarray*}
as was to be shown.
\end{proof}

By Proposition 7.8 and Lemma 7.10 we obtain
$$
\thetadot(\gamma) = - \inter\left(\int_CW\wedge
\omegaone\right)\theta(\gamma)
= \inter\left(\int_C
2\Re\left(\left(N(\omega'\omega')- 2\omega'_{(1)}\omega'_{(1)}\right)
\mudot\right)\right)\theta(\gamma).
$$
This completes the proof of Theorem 7.2.\qed

\section{The Maurer-Cartan Form on the Teichm\"uller Space}

We conclude the paper by discussing what we have shown above  
in the context of the geometry of the moduli space of compact  Riemann
surfaces. \par Let $\Sg$ be a $2$-dimensional oriented connected closed 
$C^\infty$ manifold of genus $g \geq 1$, $p_0 \in \Sg$, and 
$v_0 \in (T_\bR\Sg)_{p_0}\setminus\{0\}$. We may define 
the fundamental group $\pi_1(\Sg, p_0, v_0)$ as in \S4, 
which is also isomorphic to the free group $F_{2g}$. 
Throughout this section we fix an isomorphism and 
identify $\pi_1(\Sg, p_0, v_0)$ with $F_{2g}$, 
so that we identify $H = H_1(F_{2g}; \bR)$ with $H_1(\Sg; \bR) 
\overset\vartheta\cong H^1(\Sg; \bR)$. \par
The universal covering space $\widetilde{\moduli}$ of the space  
$\moduli$ is the Teichm\"uller space $\mathcal{T}_{g, 1}$ of triples 
of genus $g$. It is, by definition, the moduli space of 
quadruples $(C, P_0, v, [f])$ of genus $g$, where $(C, P_0, v)$ 
is a triple of genus $g$ and $[f]$ the isotopy class of 
an orientation-preserving diffeomorphism $f: (\Sg, p_0, v_0)\to 
(C, P_0, v)$. For any isomorphism class $[C, P_0, v, [f]]$ 
of quadruples 
we may consider the Magnus expansion
\begin{eqnarray*}
&&\vert f_*\vert\inv\circ \theta^{(C, P_0, v)}\circ f_*: \\
&&F_{2g} = \pi_1(\Sg, p_0, v_0) \overset{f_*}\longrightarrow 
\pi_1(C, P_0, v) \overset{\theta^{(C, P_0, v)}}\longrightarrow
1 + \T_1(H_1(C; \bR)) 
\overset{\vert f_*\vert\inv}\longrightarrow 1 + \T_1,
\end{eqnarray*}
where $\vert f_*\vert = f_*: H = H_1(\Sg) \to H_1(C)$ is the induced 
homomorphism of the diffeomorphism $f$. Consequently we obtain 
a real analytic equivariant map
\begin{equation}
\theta: \widetilde{\moduli} = \mathcal{T}_{g, 1} \longrightarrow
 \Theta_{2g}, \quad
[C, P_0, v, [f]] \mapsto 
\vert f_*\vert\inv\circ \theta^{(C, P_0, v)}\circ f_*
\label{8.1}
\end{equation}
with respect to the mapping class group. 
We call it {\bf the harmonic Magnus expansion} on the universal 
family of compact Riemann surfaces over the moduli space $\moduli$. 
\par
Here we should recall how the cotangent space $T^*_{[C, P_0, v]}
\moduli$ at a point $[C, P_0, v] \in \moduli$ is identified with 
the space of quadratic differentials $H^0(C; \mathcal{O}_C((T^*C)^{\otimes
2}\otimes [P_0]^{\otimes 2}))$. Here $[P_0]$ is the holomorphic line 
bundle on $C$ defined by the divisor $P_0$. For simplicity we write
$$
H^q(C; \lambda K + \nu P_0) :=
H^q(C; \mathcal{O}_C((T^*C)^{\otimes\lambda}\otimes [P_0]^{\otimes \nu}))
$$
for $\lambda, \nu \in \bZ$ and $q = 0, 1$. As in \S7 let 
$(C_t, P^t_0, v_t)$ be a $C^\infty$ family of triples satisfying 
$(C_t, P^t_0, v_t)\Bigr\vert_{t=0} = (C, P_0, v)$, and 
$f^t: (C, P_0, v) \to (C_t, P^t_0, v_t)$ a $C^\infty$ family of 
diffeomorphisms satisfying the condition (7.1). 
The Kodaira-Spencer isomorphism $T_{[C, P_0, v]}
\moduli = H^1(C; -K-2P_0)$ maps the tangent vector 
$\frac{d}{dt}\Bigr\vert_{t=0}(C_t, P^t_0, v_t)$ to the Dolbeault 
cohomology class $[\mudot]$. 
The Serre duality $H^0(C; 2K+2P_0) \otimes H^1(C; -K-2P_0)\to 
\bC$, $q\otimes \mudot \mapsto \int_Cq\mudot$, gives a natural 
isomorphism
\begin{equation}
T^*_{[C, P_0, v]}\moduli = H^0(C; 2K+2P_0). 
\label{8.2}
\end{equation}
\par
Let $\modg$ be the moduli space of compact Riemann surfaces $C$ 
of genus $g$, and $\family$ the moduli of pointed compact Riemann 
surfaces $(C, P_0)$. The forgetful map $\pi: \family \to \modg$, 
$[C, P_0] \mapsto [C]$, can be regarded as the universal family 
of compact Riemann surfaces over $\modg$. The space $\moduli$ is 
obtained by deleting the zero section from the relative tangent 
bundle of $\pi$, $\moduli = T_{\family/\modg} \setminus 0(\family)$. 
We denote the projection of the relative tangent bundle by $\varpi: 
\moduli \to \family$, which coincides with the forgetful map $[C, P_0, v]
\mapsto [C, P_0]$. Similarly we have natural isomorphisms
\begin{equation}
T^*_{[C]}\modg = H^0(C; 2K), \quad\text{and}\quad
T^*_{[C, P_0]}\family = H^0(C; 2K+P_0). 
\label{8.3}
\end{equation}
We may regard $H$ and $H_\bC := H\otimes \bC$ 
as flat vector bundles over the moduli space $\modg$.
It follows from
Theorem 7.2 and (2.7) 
\begin{eqnarray*}
&& \langle \theta^*\eta, \mudot\rangle = \left\langle\theta^*\eta,
\frac{d}{dt}\Bigr\vert_{t=0}(C_t, P^t_0, v_t)\right\rangle\\
&=& \,\,\inter\left(\int_C 2\Re\left(\left(N(\omega'\omega')-
2\omega'_{(1)}\omega'_{(1)}\right)\mudot\right)\right)
 \in \operatorname{Der}(\T).
\end{eqnarray*}
By (4.6) we obtain
\begin{theorem}\label{Theorem 8.1}
$$
\theta^*\eta_p = N(\omega'\omega')_{(p+2)} 
+ \overline{N(\omega'\omega')_{(p+2)}} \in 
(T^*_\bR \moduli)_{[C, P_0, v]}\otimes H^{\otimes (p+2)}
$$
for any $p \geq 1$. 
\end{theorem}
In fact, $N(\omega'\omega')$ is invariant under the action of $\varepsilon$.
\par
We may consider the pullback of the $p$-cocycle $Y_p$ by the map $\theta$
$$
\theta^*Y_p \in C^*(K_{p+1}; \Omega^*(\moduli; H^{\otimes
(p+2)})).
$$
Applying Theorem 3.2 to the map $\theta$ we have 
\begin{theorem}\label{Theorem 8.2}
$$
[\theta^*Y_p] = \frac1{(p+2)!}(-1)^{\frac12p(p+1)} m_{0, p+2} \in 
H^p(\moduli; H^{\otimes (p+2)})
$$
under the isomorphism $H^p(C^*(K_{p+1}; \Omega^*(\moduli; H^{\otimes
(p+2)}))) = H^p(\moduli; H^{\otimes(p+2)})$.
\end{theorem}
This can be interpreted the $p$-cocyle $\theta^*Y_p$ is 
a canonical combinatorial family of differential forms 
representing the twisted Morita-Mumford class $m_{0, p+2}$. 
From \cite{Mo3} and \cite{KM1} the $i$-th Morita-Mumford class $e_i$ is 
obtained by contracting the coefficients of $m_{0, 2i+2}$ 
using the intersection form of the surface.\par

In the case $p \neq 2$, as was shown in Lemma 7.1, 
$N(\omega'\omega')_{(p+2)}$ has a pole of order $\leq 1$ 
at the point $P_0$. Hence $\theta^*\eta_p$ can be regarded as 
a twisted $1$-form on $\family$, $\theta^*\eta_p \in 
\Omega^1(\family; H^{\otimes (p+2)})$ for any $p \neq 2$. \par

Finally we discuss more about the twisted closed $1$-form 
$\theta^*\eta_1 \in \Omega^1(\family; H^{\otimes 3})$. 
For the rest of the paper we suppose $g \geq 2$. 
Let $\Mg$ and $\Mgstar$ be the mapping class groups of the surface 
$\Sg$ and the pointed surface $(\Sg, p_0)$, respectively. 
The space $\modg$ is the quotient of the Teichm\"uller space 
of genus $g$ by the natural action 
of the group $\Mg$. Hence, for any $\bR[\Mg]$-module 
$M$, we have a natural isomorphism $H^*(\modg; M) = 
H^*(\Mg; M)$. Similarly we have $H^*(\family; M) = 
H^*(\Mgstar; M)$. \par
Morita introduced the extended Johnson homomorphisms 
$\tilde k \in H^1(\Mgstar; \La^3H)$ and $\tilde k \in 
H^1(\Mg; \La^3H/H)$ in \cite{Mo2}. Here we regard $H$ as 
an $Sp_{2g}(\bR)$-submodule of $\La^3H$ through the injection
$\frak q^H: Z \in H \mapsto Z\wedge I \in \La^3H$ given by the wedge product 
by the intersection form $I$.  
The forgetful homomorphism $\varpi^*: H^1(\Mgstar; 
\La^3H) \to H^1(\Mgone; \La^3H)$ is an isomorphism, 
and maps $\tilde k$ to the cohomology class of the 
first Johnson map $[\tau^\theta_1] = h_1$.\par

As in \S4 we choose a symplectic basis 
$\{X_i, X_{g+i}\}^g_{i=1}$ of $H = H_1(\Sg; \bR)$. 
For any $p \geq 2$ we regard $\La^pH$ as a 
$Sp_{2g}(\bR)$-submodule of $H^{\otimes p}$ by the injection 
\begin{equation}
\La^pH \to H^{\otimes p}, \quad
Z_1\wedge\cdots\wedge Z_p \mapsto 
{\sum}_{\sigma \in \frak S_p}(\operatorname{sgn}\sigma)Z_{\sigma(1)}
\cdots Z_{\sigma(p)}.
\label{8.4}
\end{equation}
Here it should be remarked the wedge product $\varphi_1\wedge\varphi_2$ 
of two $\T$-valued $1$-forms $\varphi_1$ and $\varphi_2$ on a Riemann 
surface does
{\bf not}  mean the wedge product of the coefficients, but only that of 
differential forms as in the preceding sections. For example we mean 
$(Z_1 dz)\wedge (Z_2d\zbar) = (Z_1Z_2)dz\wedge d\zbar =
(Z_1\otimes Z_2)dz\wedge d\zbar$ for $Z_1$ and 
$Z_2 \in H$. We define the contraction map $\frak c: \La^pH \to 
\La^{p-2}H$ by the restriction of the map $H^{\otimes p} 
\to H^{\otimes (p-2)}$, $Z_1Z_2Z_3\cdots Z_p \mapsto (Z_1\cdot Z_2)Z_3
\cdots Z_p$. For example the intersection form $I = \sum X_iX_{g+i} 
- X_{g+i}X_i = \sum X_i\wedge X_{g+i}$ satisfies $\frak c(I) = 2g$. 
Moreover we have $\frak c(Z_1\wedge Z_2\wedge Z_3) = 2((Z_1\cdot Z_2)Z_3 
+ (Z_2\cdot Z_3)Z_1 + (Z_3\cdot Z_1)Z_2)$ for any $Z_1, Z_2$, and $Z_3 
\in H$. \par

Now we may regard $\theta^*\eta_1 \in \Omega^1(\family; \La^3H)$. 
In fact, the $(1, 0)$-part of $\theta^*\eta_1$ is 
\begin{eqnarray*}
N(\omega'\omega')_{(3)} &=& 2N(\omega'_{(1)}\omega'_{(2)}) 
= -2\sqrt{-1}N(\omega'_{(1)}\dd\Phi(\omegaone\wedge\omegaone)) \\
&=& -\sqrt{-1}N(1- (23))(\omega'_{(1)}\dd\Phi(\omegaone\wedge\omegaone))
\\ 
&=& -\sqrt{-1}({\sum}_{\sigma\in \frak
S_3}(\operatorname{sgn}\sigma)\sigma)
(\omega'_{(1)}\dd\Phi(\omegaone\wedge\omegaone)), 
\end{eqnarray*}
which has coefficients in $\La^3H\otimes \bC$. 
From (2.6) abd the results of Morita 
quoted above we obtain
\begin{equation}
[\theta^*\eta_1] = -\tilde k \in H^1(\family; \La^3H). 
\label{8.5}
\end{equation}
Following \cite{Mo3} we write $U := \La^3H/H = \Coker\,\frak q^H$. If we
denote by
$\frak p^U: \La^3H \to U$ the natural projection, 
$\frak p^H:= \frac{1}{2g-2}\frak c: \La^3H \to H$, and
$\frak q^U: = 1 -  \frak q^H\frak p^H : U \to \La^3H$, 
then we have $\frak p^H\frak q^H = 1_H$ and these maps give an 
$Sp_{2g}(\bR)$-equivariant decomposition
\begin{equation}
\La^3H = H \oplus U.
\label{8.6}
\end{equation}
As was shown in \cite{Mo2} the forgetful map $H^1(\Mg; U) \to H^1(\Mgstar;
U)$  is injective and maps $\tilde k$ to $\frak p^U\tilde k$. 
We define 
$$
\eta^H_1 := \frak p^H\theta^*\eta_1 \in \Omega^1(\family; H) \quad
\text{and}\quad
\eta^U_1 := \frak p^U\theta^*\eta_1 \in \Omega^1(\family; U).
$$
The closed $1$-form $\eta^U_1$ represents the image of the extended 
Johnson homomorphism on $\Mg$ by the forgetful map. 
\begin{theorem}\label{Theorem 8.3} For any pointed Riemann surface $(C,
P_0)$  the value of the $(1, 0)$-part of $\eta^U_1$ at $[C, P_0] \in
\family$, 
$\frak p^U N(\omega'\omega')_{(3)} \in H^0(C; 2K+P_0)\otimes U$, 
is smooth at $P_0$, and independent of the choice of the point $P_0$. 
In other words, $\eta^U_1$ can be regarded as a $1$-form on the moduli 
space $\modg$, and represents the extended Johnson homomorphism $-\tilde
k$  on the mapping class group $\Mg$.
\end{theorem}

To prove the theorem we fix a pointed Riemann surface $(C, P_0)$, 
and choose the basis $\{\xi_i\}^{2g}_{i=1}$ 
of the real harmonic $1$-forms on $C$, whose cohomology classes 
form the dual basis of $\{X_i\}^{2g}_{i=1}$ as in \S5. 
We have $\omegaone = \sum^g_{i=1}(\xi_iX_i + \xi_{g+i}X_{g+i})$.
We denote 
$$
B = B_C := \frac1{2g}\omegaone\cdot\omegaone 
= \frac1g {\sum}^g_{i=1}\xi_i\wedge\xi_{g+i} \in \Omega^2(C).
$$
Clearly $\int_CB = 1$. Moreover $B$ is a volume form on $C$. 
In fact, if we choose a basis $\{\psi_i\}^g_{i=1}$ of the space 
$H^0(C; K)$ satisfying the condition $\int_C\psi_i\wedge
\overline{\psi_j} = \delta_{i,j}$, $1\leq i, j \leq g$, then 
we have $B = \frac{\sqrt{-1}}{2g}{\sum}^g_{i=1}\psi_i\wedge
\overline{\psi_i}$ because $\{\psi_i, \frac{\sqrt{-1}}2
\overline{\psi_i}\}^g_{i=1}$ is a symplecic basis of the space 
$H\otimes \bC$. Since the complete linear system $\vert K\vert$ 
of the canonical divisor on the Riemann surface $C$ has no basepoint, $B$
does not vanish  at any point $P \in C$. \par
We denote by 
$\Omega^2_0(C)$ the kernel of the map $\int_C: \Omega^2(C)\otimes \bC
 \to \bC$. 
The key to the proof of the theorem is the fact the Green operator 
$\Phi: A^2(C) \to A^0(C)/\bC$ restricted to the kernel $\Omega^2_0(C)$ 
has values in the smooth functions $(\Omega^0(C)\otimes \bC)/\bC$ and 
is independent of the choice of the point $P_0$, which we denote
$$
\Phi_0 = \Phi^C_0: \Omega^2_0(C) \to (\Omega^0(C)\otimes \bC)/\bC.
$$
For example, since $\int_C (\omegaone\wedge\omegaone -BI) = I - I 
= 0$, the $1$-form $\ast d\Phi_0(\omegaone\wedge\omegaone -BI) = \ast
d\Phi(\omegaone\wedge\omegaone -BI)$  is smooth at $P_0$, 
and independent of the point $P_0$.
We denote the $(1, 0)$-part of $\eta^H_1$ by 
$$
q = q^{(C, P_0)} := \frak p^HN(\omega'\omega')_{(3)}
= 2\frak p^HN(\omega'_{(1)}\omega'_{(2)}).
$$

\begin{lemma}\label{Lemma 8.4} The covariant tensor $q^{(C, P_0)} -
2\omega'_{(1)}\ast\dd\Phi B
\in C^\infty(C \setminus\{P_0\}; (T^*C)^{\otimes 2})\otimes H$ is smooth at
$P_0$,  and independent of the choice of $P_0$.
\end{lemma}
\begin{proof} One computes
\begin{eqnarray*}
&& \frak c N(\omega'_{(1)}\omega'_{(2)}) \\
&=& 2g\omega'_{(1)}\ast\dd\Phi B 
+ 2 \sum\xi'_i\ast\dd\Phi(\xi_{g+i}\wedge\omegaone) 
- 2 \sum\xi'_{g+i}\ast\dd\Phi(\xi_{i}\wedge\omegaone) \\
&=& 2g\omega'_{(1)}\ast\dd\Phi B \\
&&+ 2 \sum\xi'_i\ast\dd\Phi_0(\xi_{g+i}\wedge(\omegaone-\xi_iX_i)) 
- 2 \sum\xi'_{g+i}\ast\dd\Phi_0(\xi_{i}\wedge(\omegaone-\xi_{g+i}X_{g+i})) \\
&&- 2
\sum(\xi'_iX_i+\xi'_{g+i}X_{g+i})\ast\dd\Phi(\xi_{i}\wedge\xi_{g+i})\\
&=& (2g-2)\omega'_{(1)}\ast\dd\Phi B \\
&&+ 2 \sum\xi'_i\ast\dd\Phi_0(\xi_{g+i}\wedge(\omegaone-\xi_iX_i)) 
- 2 \sum\xi'_{g+i}\ast\dd\Phi_0(\xi_{i}\wedge(\omegaone-\xi_{g+i}X_{g+i})) \\
&&- 2
\sum(\xi'_iX_i+\xi'_{g+i}X_{g+i})\ast\dd\Phi_0(\xi_{i}\wedge\xi_{g+i}-B).\\
\end{eqnarray*}
This means $\frak p^HN(\omega'\omega')_{(3)} - 2\omega'_{(1)}\ast\dd\Phi B$ 
is smooth at $P_0$, 
and independent of the choice of $P_0$.\qed
\end{proof}

\begin{proof}[Proof of Theorem 8.3] 
For any $Z \in H$ we have $Z\wedge I =
N(ZI)$.  The $(1, 0)$-part of $\eta^U_1$,
$N(\omega'\omega')_{(3)} - q^{(C, P_0)}\wedge I$, is congruent to 
$2N(\omega'_{(1)}\omega'_{(2)}) - 2 \omega'_{(1)}\wedge 
\ast\dd\Phi BI$ modulo the tensors on $C$ smooth at $P_0$, and independent of
the choice of $P_0$. Now the difference 
$N(\omega'_{(1)}\omega'_{(2)}) - \omega'_{(1)}\wedge 
\ast\dd\Phi BI
= N(\omega'_{(1)}\ast\dd\Phi_0(\omegaone\wedge\omegaone - BI))$
is smooth at $P_0$ and independent of the choice of $P_0$. 
This proves the theorem.
\end{proof}

In \cite{H} Harris introduced the harmonic volume $I_C$ of a compact
Riemann surface $C$, and computed the first variation of $I_C$. 
The $1$-form $\eta^U_1$ coincides with the first variation. 
As was pointed out in \cite{H}, $\eta^U_1$ vanishes along the
hyperelliptic  locus $\mathcal{H}_g \subset \modg$. One can see this fact
by considering the action of the hyperelliptic involution. This implies
all the differential  forms representing the Morita-Mumford classes
derived from 
 $\eta^U_1$ vanish along the locus $\mathcal{H}_g$.
The harmonic volume of a
hyperelliptic Riemann surface, however, is non-trivial, 
as was pointed out by Tadokoro \cite{T1}. 
It is a locally constant function on the locus and 
computed completely in \cite{T1}. See also \cite{T2}.\par

In view of a theorem of Morita \cite{Mo3} there exist unique equivariant
linear maps $\alpha_0$ and 
$\alpha_1 \in \Hom((\La^3H)^{\otimes 2}, \bR)^{Sp_{2g}(\bR)}$ such that 
$e^J := {\alpha_0}_*({\theta^*\eta_1}^{\otimes 2})$ and ${e_1}^J :=
{\alpha_1}_*({\theta^*\eta_1}^{\otimes 2}) \in \Omega^2(\family)$
represent the first Chern class of $T_{\family/\mathbb{M}_g}$ and 
the first Morita-Mumford class, respectively. 
One can construct the form $e^J_1$ from the $1$-form 
$\eta^U_1$, and so regard it as a $2$-form
on $\modg$.  Hain and Reed \cite{HR} 
introduced the same $1$-form $\eta^U_1$ and studied the difference $e^J_1
- 12c_1(\lambda, L^2)$ in detail, where $\lambda$ is the Hodge line bundle
over the Siegel upper half space $\frak H_g$.\par 
The Chern form $e^J$ seems to be related to Arakelov's admissible
metric.  Let $B$  be the
volume form on a compact Riemann surface $C$ introduced above, and 
$h$ the function on $\family\times_{\mathbb{M}_g}\family$ 
satisfying the conditions $\frac1{2\pi\sqrt{-1}}\dd\dc
h\bigr\vert_{C\times\{P_0\}} = B -
\delta_{P_0}$ and 
$\int_C\left(h\bigr\vert_{C\times\{P_0\}}\right)B = 0$. 
Then we have 
$$
\left(\frac1{2\pi\sqrt{-1}}\dd\dc h\right)\biggr\vert_{\text{diagonal}} 
= e^J + \frac1{(2-2g)^2}({e_1}^F - {e_1}^J
) \,\,\in\,\, \Omega^2(\family).
$$
Here ${e_1}^F := \fiber (e^J)^2 \in \Omega^2(\mathbb{M}_g)$. 
The proof will be given in \cite{Kaw4}.

\bibliographystyle{amsplain}

\end{document}